\documentclass[11pt,preprint]{elsarticle}

\makeatletter
\def\ps@pprintTitle{%
 \let\@oddhead\@empty
 \let\@evenhead\@empty
 \def\@oddfoot{\centerline{\thepage}}%
 \let\@evenfoot\@oddfoot}
\makeatother

\usepackage{microtype}
\usepackage[]{algorithm2e}
\usepackage{chngpage}
\usepackage[english]{babel}
\usepackage[T1]{fontenc}
\usepackage[utf8]{inputenc}

\usepackage{amsfonts}
\usepackage{amsmath}
\usepackage{amssymb}
\usepackage{mathtools}
\usepackage[binary-units=true]{siunitx}
\usepackage{breqn}
\usepackage{stmaryrd}
\usepackage{sansmath}
\usepackage{amsthm}
\usepackage{bm}

\usepackage{caption}
\usepackage{subcaption}
\usepackage{diagbox}

\usepackage{rotating}

\usepackage{geometry}
\newgeometry{vmargin={25mm}, hmargin={25mm,25mm}}

\usepackage{multicol}
\usepackage{booktabs}
\usepackage{lscape}
\makeatletter
\newcommand\tabfill[1]{%
\dimen@\linewidth%
\advance\dimen@\@totalleftmargin%
\advance\dimen@-\dimen\@curtab%
\parbox[t]\dimen@{#1\ifhmode\strut\fi}%
}

\usepackage{multirow}

\usepackage{dcolumn}
\newcolumntype{d}[1]{D{.}{.}{#1}}

\usepackage{tikz}
\usepackage[]{graphicx}

\usetikzlibrary{shapes.geometric}

\newcommand*{\ldcbrace}{\{\mskip-5mu\{}
\newcommand*{\rdcbrace}{\}\mskip-5mu\}}

\renewcommand{\@algocf@capt@plain}{above}

\definecolor{clr1}{rgb}{0.031, 0.270, 0.580}
\definecolor{clr2}{rgb}{0.129, 0.443, 0.709}
\definecolor{clr3}{rgb}{0.258, 0.572, 0.776}
\definecolor{clr4}{rgb}{0.419, 0.682, 0.839}
\definecolor{clr5}{rgb}{0.619, 0.792, 0.882}

\usepackage[unicode=true,
            bookmarks=true,
            bookmarksnumbered=false,
            bookmarksopen=false,
            breaklinks=true,
            pdfborder={0 0 1},
            backref=false,
            colorlinks=true,
            linkcolor=black,
            urlcolor=theblue,
            citecolor=theblue,
            hyperfootnotes=false]{hyperref}
\usepackage[mathlines]{lineno}
\hypersetup{pdftitle={},
            pdfauthor={The Authors}}

\let\oldref\ref
\renewcommand{\ref}[1]{(\oldref{#1})}

\newtheorem{theorem}{Theorem}[section]

\newtheorem{lemma}[theorem]{Lemma}

\allowdisplaybreaks

\usepackage{xcolor}

\newcommand{\colornameedits}{black}

\usepackage{placeins}
\DeclareMathOperator\arctanh{arctanh}
\begin{document}
\newcommand{\hwplotFR}{\raisebox{2pt}{\tikz{\draw[black,solid,line width=0.9pt](0,0) -- (5mm,0);}}}
\newcommand{\hwplotHFR}{\raisebox{2pt}{\tikz{\draw[black,dashdotted,line width=0.9pt](0,0) -- (5mm,0);}}}
\newcommand{\hwplotEFR}{\raisebox{2pt}{\tikz{\draw[black,dotted,line width=0.9pt](0,0) -- (5mm,0);}}}

\newcommand{\hwplotDG}{\raisebox{2pt}{\tikz{\draw[black,solid,line width=0.9pt](0,0) -- (5mm,0);}}}
\newcommand{\hwplotSD}{\raisebox{2pt}{\tikz{\draw[black,dashed,line width=0.9pt](0,0) -- (5mm,0);}}}
\newcommand{\hwplotHU}{\raisebox{2pt}{\tikz{\draw[black,dotted,line width=0.9pt](0,0) -- (5mm,0);}}}

\newcommand{\hwplotDGmarker}{\raisebox{0pt}{\tikz{\node[circle,draw,fill=none,text width=2mm, inner sep=0pt] at (0,0){};}}}
\newcommand{\hwplotSDmarker}{\raisebox{0pt}{\tikz{\node[isosceles triangle,draw,fill=none,rotate=90,text width=2mm, inner sep=0pt] at (0,0){};}}}
\newcommand{\hwplotHUmarker}{\raisebox{0pt}{\tikz{\node[regular polygon,regular polygon sides=4,draw,fill=none,text width=1.25mm, inner sep=0pt] at (0,0){};}}}

\newcommand{\hwplotDGmarkerfilled}{\raisebox{0pt}{\tikz{\node[circle,draw,fill=black,text width=2mm, inner sep=0pt] at (0,0){};}}}
\newcommand{\hwplotSDmarkerfilled}{\raisebox{0pt}{\tikz{\node[isosceles triangle,draw,fill=black,rotate=90,text width=2mm, inner sep=0pt] at (0,0){};}}}
\newcommand{\hwplotHUmarkerfilled}{\raisebox{0pt}{\tikz{\node[regular polygon,regular polygon sides=4,draw,fill=black,text width=1.25mm, inner sep=0pt] at (0,0){};}}}

\title{Hybridized Formulations of Flux Reconstruction Schemes for Advection-Diffusion Problems}

\begin{frontmatter}
\author[one]{Carlos A. Pereira\corref{mycorrespondingauthor}}
\cortext[mycorrespondingauthor]{Corresponding author}
\ead{carlos.pereira@concordia.ca}

\author[one]{Brian C. Vermeire}
\ead{brian.vermeire@concordia.ca}

\address[one]{Department of Mechanical, Industrial, and Aerospace Engineering. \\Concordia University. \\Montreal, QC. Canada}

\begin{abstract}
 We present the hybridization of flux reconstruction methods for advection-diffusion problems. Hybridization introduces a new variable into the problem so that it can be reduced via static condensation. This allows the solution of implicit discretizations to be done more efficiently. We derive an energy statement from a stability analysis considering a range of correction functions on hybridized and embedded flux reconstruction schemes. Then, we establish connections to standard formulations. We devise a post-processing scheme that leverages existing flux reconstruction operators to enhance accuracy for diffusion-dominated problems. Results show that the implicit convergence of these methods for advection-diffusion problems can result in performance benefits of over an order of magnitude. In addition, we observe that the superconvergence property of hybridized methods can be extended to the family of FR schemes for a range of correction functions.
\end{abstract}
\end{frontmatter}
\section{Introduction}
Flux reconstruction (FR) is a family of high-order numerical methods that can recover existing high-order schemes, including discontinuous Galerkin (DG)~\cite{cockburnTBVRungeKuttaLocal1989,cockburnTVBRungeKuttaLocal1989}, spectral difference (SD)~\cite{liu2006discontinuous} and many others via the choice of a correction function. A range of energy-stable correction functions was determined by Vincent et al.~\cite{vincentNewClassHighOrder2011} for tensor-product formulations by a continuous scalar parameter. These include the scheme defined by Huynh~\cite{huynhFluxReconstructionApproach2007} in connection with the DG-SEM method, generally known as FR$_{g_2}$ or FR$_{HU}$-method. These correction functions provide penalization of flux discontinuities at interfaces between elements. Linear stability proofs of FR methods have been obtained for linear advection~\cite{vincentNewClassHighOrder2011,castonguayEnergyStableFlux2013,pereiraSpectralPropertiesHighOrder2020} and linear diffusion problems~\cite{quaegebeurStabilityEnergyStable2019,quaegebeurStabilityEnergyStable2020} with different types of viscous Riemann solvers. These include formulations with local-discontinuous Galerkin (LDG)~\cite{cockburnLocalDiscontinuousGalerkin1998a}, Bassi and Rebay II (BR2)~\cite{bassi1997high}, and interior penalty (IP)~\cite{arnoldInteriorPenaltyFinite1982} formulations. 

While FR methods have been generally used with explicit time stepping, stiff problems benefit greatly from implicit methods due to their less strict stability constraints. However, they can be significantly more expensive per step as they require the solution of prohibitively large systems of nonlinear equations, especially for high-order discretizations. Cockburn et al.~\cite{cockburnUnifiedHybridizationDiscontinuous2009} presented the hybridization of discontinuous Galerkin (HDG) methods and its connection to other finite-element approaches such as the Raviart-Thomas (RT) and Brezzi-Douglas-Merini (BDM) methods. These introduce an additional unknown into the problem, the trace variable, which acts as a boundary communicator between adjacent elements. Hence, the conservation law is discretized into a set of local problems, and only a system involving this trace is to be solved after static condensation. The behaviour of these schemes depends on the function space where this trace variable lives, as well as the choice of the so-called stabilization parameter. In terms of function spaces, these include methods with discontinuous trace polynomials, generally referred to as HDG, and methods where this trace is set to be continuous on the skeleton of the computational domain, known as the embedded DG (EDG) methods~\cite{cockburnAnalysisEmbeddedDiscontinuous2009,nguyenClassEmbeddedDiscontinuous2015}. Hybridized high-order methods effectively reduce the scaling of the implicit system, generally $(p+1)^{2d}$ to a lower dimension, i.e  $(p+1)^{2(d-1)}$, where $p$ is the polynomial degree of the solution and $d$ is the number of dimensions.

For pure linear advection problems with discontinuous traces, a direct connection to DG methods has been established in the literature for upwind-like stabilization~\cite{fernandezEntropystableHybridizedDiscontinuous2019a,nguyenImplicitHighorderHybridizable2009c}. For diffusion problems, it has been shown that no finite stabilization can recover existing conventional discretizations~\cite{cockburnUnifiedHybridizationDiscontinuous2009}. One well-known property of HDG methods is that they possess optimal $p+1$ convergence of the solution and flux variables when diffusion operators dominate and that the solution has a superconvergent behaviour with respect to a projection~\cite{cockburnConditionsSuperconvergenceHDG2012,cockburnSuperconvergentDiscontinuousGalerkin2009a}. This allows application of a post-processing scheme to enhance the order of accuracy of the solution from the conventional $p+1$ to $p+2$~\cite{cockburnSuperconvergentDiscontinuousGalerkin2009a,cockburnSuperconvergentLDGhybridizableGalerkin2008}. Hybridized methods have been shown to be applicable to a wide variety of problems. These include linear convection~\cite{jaustFESTUNGMATLABGNU2018} and convection-dominated~\cite{nguyenImplicitHighorderHybridizable2009,nguyenImplicitHighorderHybridizable2009c} problems, as well as nonlinear problems pertinent to computational fluid dynamics, such as incompressible~\cite{giacominiTutorialHybridizableDiscontinuous2020a,cesmeliogluAnalysisHybridizableDiscontinuous2017} and compressible flows~\cite{vila-perezHybridisableDiscontinuousGalerkin2021}, and turbulent flows~\cite{fernandezImplicitLargeeddySimulation2016}. 

While these methods have been shown to display superior performance and accuracy,  hybridization has not been explored in the context of flux reconstruction schemes until our recent publication~\cite{pereiraPerformanceAccuracyHybridized2022}. There, it was shown that for advection problems, significant speedup factors can be obtained with similar numerical error behaviour to conventional FR methods for a variety of correction functions. However, to the authors' knowledge, their properties in the context of convection-diffusion problems have not been studied. In this paper, we explore the stability, accuracy and performance of hybridized flux reconstruction schemes. Specifically, we consider the hybridized flux reconstruction (HFR) and embedded flux reconstruction (EFR) methods and compare these properties against conventional implicit FR-LDG discretizations. 

This manuscript is structured as follows. In Section 1, we present the hybridization framework and implementation of advection-diffusion problems. Then, we analyze the stability mechanisms of hybridized FR methods via an energy analysis and establish the connection to standard FR formulations in Section 2. In Section 3, we present a post-processing scheme that leverages existing operators from FR methods and is applicable to all schemes within this family. We then analyze these properties via linear and nonlinear numerical examples, including steady-state and unsteady problems. We finalize with conclusions and recommendations for future work.

\section{The Hybridized Flux Reconstruction Method}\label{sec:hfr}

\subsection{Preliminaries}
Consider the two-point boundary value problem  
\begin{equation}
 \frac{\partial u}{\partial t} + \nabla \cdot \bm F(u, \nabla u) = 0~~\text{in}~\Omega,
\end{equation}
where $\Omega$ is a bounded subset of $\mathbb{R}^d$ with boundary $\partial\Omega\in \mathbb{R}^{d-1}$ and $d$ dimensions, $u$ is the conserved quantity, $\bm F=\bm F(u, \nabla u)$ is the flux and $t$ is time. To discretize this problem, we rewrite it as a system of first-order ordinary differential equations
\begin{subequations}
\begin{align}
  \frac{\partial u}{\partial t} + \nabla \cdot \bm F(u, \bm q) &= 0,\\
  \bm q - \nabla u &=0,
\end{align}\label{eq:systempde}\end{subequations}
where $\bm q$ is an auxiliary variable referring to the gradient of the conserved variable. 

We maintain the notation used in our previous work~\cite{pereiraPerformanceAccuracyHybridized2022}, where we define $\mathcal{T}_h$ to be the partition of $\Omega$ into $N$ nonoverlapping, conforming elements $\Omega_k$, each with boundary $\partial\Omega_k=\{f\}$. Define also $\partial\mathcal{T}_h$ to be the collection of all element borders such that $\partial \mathcal{T}_h=\{\partial \Omega_k: \Omega_k \in \mathcal{T}_h\}$. Here, every face is counted from the point of view of every element. \textcolor{\colornameedits}{Furthermore, consider the collection of all unique faces $\varepsilon^h=\varepsilon^h_\partial\cup\varepsilon^h_0=\{\bar f\}$ to be the union of all unique boundary $(\varepsilon^h_\partial)$ and interior $(\varepsilon^h_0)$ faces in the computational domain. Here, each face is counted only once}. Hence, two interior faces from $\partial\mathcal{T}_h$ have a single corresponding face or vice-versa $\bar f\in \varepsilon^h,~\bar f = (f\in\partial\Omega_k) \cap \varepsilon^h$.

The first step in the implementation is to map each $\Omega_k$ to a reference element $\tilde\Omega$. Then, we can make use of invertible one-to-one mapping functions $\bm {\mathcal{M}}_k( \tilde{\bm x})$ to convert quantities between physical and reference space within each element. These are obtained via
\begin{equation}
  \bm x = \bm {\mathcal{M}}_k( \tilde{\bm x}) = \sum_{i=1}^{N_g} M_i( \tilde{\bm x}) \bm x^g_i,
\end{equation}
where ${\bm x}$ is the physical coordinate of a given point in $\Omega_k$, and $M_i$ is a shape function associated with one of the $N_g$ mapping points $\{\bm x^g_i\}$. Define the Jacobian matrix of these transformations by $\bm J_k(\tilde{\bm x})$ and its determinant by $J_k(\tilde{\bm x})$. \textcolor{\colornameedits}{These geometric parameters allow us to rewrite the conservation law considering the approximated physical solution $u^h=\cup_{k=1}^N u^h_k$ and flux $\bm F^h=\cup_{k=1}^N\bm F^h_k$ in reference space such that for a time-invariant formulation~\cite{zwanenburgEquivalenceEnergyStable2016a}}
\begin{align}
  \tilde{u}^h_k &= \tilde{u}^h_k(\tilde{\bm  x}, t) = J_k u^h_k(\bm{\mathcal{M}}_k(\tilde{\bm  x}), t),\\
  \tilde{\bm F}^h_k &= \tilde{\bm F}^h_k(\tilde{\bm  x}, t) = J_k \bm J_k^{-1} \bm F^h_k(\bm{\mathcal{M}}_k(\tilde{\bm  x}), t),\\
  \tilde{\bm q}^h_k &= \tilde{\bm q}^h_k(\tilde{\bm  x}, t) = \bm J_k^{T} \bm q^h_k(\bm{\mathcal{M}}_k(\tilde{\bm  x}), t), \label{eq:d9sfxda}
\end{align}
so that the evolution of the physical solution within each element satisfies
\begin{align}
  \frac{\partial u^h_k}{\partial t} + \frac{1}{J_k}\tilde \nabla \cdot\tilde{\bm F}^h_k &= 0,\\
  \tilde{\bm q}^h_k - \tilde \nabla u^h_k &=0,
  \label{eq:3289u4}
\end{align}
where $\tilde\nabla$ is the divergence operator in reference space. 
\subsection{Implementation}
\begin{figure}[ht]
  \centering
  \includegraphics[width=0.6\textwidth]{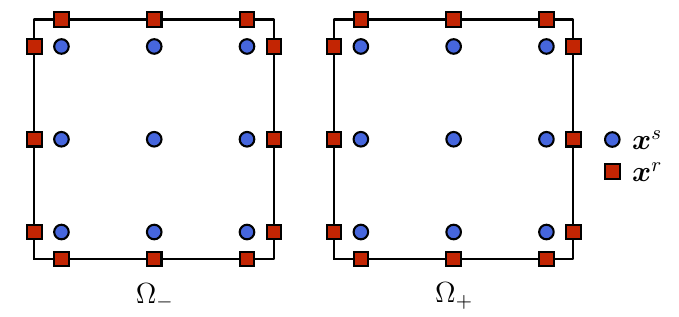}
  \caption[]{Solution and flux point locations for two neighbouring quadrilateral elements for a $p=2$ FR method}
  \label{fig:fluxpointsfr}
\end{figure}
To obtain a discretization of order $p+1$, we place $N_s$ solution points $\{\tilde{\bm x}_i^s\}_{i=1}^{N_s}$ inside each element and $N_r$ flux points $\{\tilde{\bm x}_{f,i}^r\}_{i=1}^{N_r}$ on each of the $N_f$ element faces. At each face, we define functions such that outward unit normal vectors can be obtained via $\tilde{\bm n}_f^m = \tilde{\bm n}_f(\tilde{\bm x}_{f,m}^r)$ and ${\bm n}_{k,f}^m = {\bm n}_{k,f}(\tilde{\bm x}_{f,m}^r)$ in reference and physical space, respectively. 

The conserved variable can be represented within a given element via interpolation with nodal basis functions $\{\varphi_i(\tilde{\bm x})\}_{i=1}^{N_s}$
\begin{equation}
  u^h_k (\tilde{\bm x},t) = \sum_{i=1}^{N_s} \textcolor{\colornameedits}{U_{k,i}}(t) \varphi_i(\tilde{\bm x}),
  \label{eq:frsolution}
\end{equation}
resulting in degree-$p$ discontinuous polynomials. Note that a discontinuous flux $\tilde{\bm F}^{hD}_{k}$ can also be represented using the same basis functions. In addition to the conserved variable, introduce an approximation to $u^h$ on the skeleton of the computational grid such that at a face $\bar f$, a degree $p$  polynomial can be obtained via
\begin{equation}
  \hat{u}^h_{\bar f} (\tilde{\bm x},t) = \sum_{i=1}^{N_r}  {\hat{U}}_{\bar f,i}(t) \phi_i(\tilde{x}),
  \label{eq:rwepsdf}
\end{equation}
which is the so-called trace variable. The hat symbol is used for the hybrid quantities throughout this work. For second-order problems, we first compute an approximation to the gradient, which we denote $\bm q^h$. This is obtained by taking the gradient of a globally continuous scalar variable. In the FR method, this is carried out by performing a reconstruction procedure which penalizes the interface jumps via correction functions. In other words, we add the following correction term to the discontinuous solution from Equation~\eqref{eq:frsolution}
\textcolor{\colornameedits}{\begin{equation}
  u_k^{hC} = \sum_{f=1}^{N_f} \sum_{m=1}^{N_r}  \tilde{\bm n}_f^m \cdot \bm g_f^m(\tilde{\bm x}) \left[\mathfrak{U}_{k,f} - u_{k,f}^h\right]_{\tilde{\bm x} = \tilde{\bm x}_{f,m}^r},
  \label{eq:solcorrection}
\end{equation}}
\textcolor{\colornameedits}{where $u_{k,f}$ is the element solution evaluated at face $f$}. The correction functions $\bm g_f^m(\tilde{\bm x})$ satisfy
\begin{equation}
  \tilde{\bm n}_f^m \cdot \bm{g}^n_l(\tilde{\bm x}_f^m) = \delta_{fl}\delta_{mn}.
  \label{eq:dsklnf}
\end{equation}
Here, subscripts indicate the face number, superscripts the node number within the face, and $\delta$ is the Kronecker delta. A brief discussion of these correction functions is given later in this work. Equation~\eqref{eq:solcorrection} requires a common value of the solution $\mathfrak{U}_{k,f}$. This quantity is typically explicitly computed via well-known methods, including LDG, BR2 and IP, that rely on weighting and directionality parameters of neighboring values of the solution. In the context of hybridized methods, this common value is computed implicitly and is simply set to $\mathfrak{U}_{k,f} = \hat{u}^h_{\bar f}$.  Hence, the auxiliary variable is a vector polynomial of degree $p$ that results from taking the gradient of the solution and its correction. For each element, we can write it as
\begin{equation}
  \tilde{\bm q}_k^h = \sum_{i=1}^{N_s} \textcolor{\colornameedits}{U_{k,i}}(t) \tilde \nabla \varphi_i(\tilde{\bm x}) - \sum_{f=1}^{N_f} \sum_{m=1}^{N_r} \tilde{\bm n}_f^m\cdot \tilde \nabla \cdot \bm g_f^m(\tilde{\bm x})\left[\hat u^h_{\bar f} - u_{k,f}^h\right]_{\tilde{\bm x} = \tilde{\bm x}_{f,m}^r} = 0,
\end{equation}
and is then mapped to physical space via Equation \eqref{eq:d9sfxda}. Following a similar procedure to computing the corrected gradient, we construct a correction term to the discontinuous flux using the same correction functions by
\begin{equation}
\tilde {\bm F}_{k}^{hC}(\tilde{\bm x}, t) = \sum_{f=1}^{N_f}\sum_{m=1}^{N_r} \bm g_f^m(\tilde{\bm x}) \left[{\tilde{ H}}(\tilde{\bm x})_{k,f}\right]_{\tilde{\bm x} = \tilde{\bm x}_{f,m}^r},
\label{eq:fluxcorrection}
\end{equation}
where the normal jump of the flux at the face is defined as follows
\begin{equation}
  \tilde{H}_{k,f}(\tilde{\bm x}) = \tilde{\hat{\bm{\mathfrak{F}}}}_{k,f}\cdot \tilde{\bm n}_f  - \tilde{\bm F}_{k,f}^{hD}\cdot\tilde{\bm n}_f.
\end{equation}
$\tilde{\bm F}_{k,f}^{hD}$ is the transformed discontinuous flux polynomial interpolated to face $f$, and the relationship between physical and reference space for the common flux is\textcolor{\colornameedits}{~\cite{zwanenburgEquivalenceEnergyStable2016a}}
\begin{equation}
  \tilde{\hat{\bm{\mathfrak{F}}}}_{k,f}\cdot \tilde{\bm n}_f = J_{k,f} \hat{\bm{\mathfrak{F}}}_{k,f}\cdot {\bm n}_{k,f}.
  \label{eq:8923i}
\end{equation}
Note that there is no constraint on applying different correction functions for the gradient and the flux in Equations~\eqref{eq:solcorrection} and~\eqref{eq:fluxcorrection}, but we choose to make use of the same in this work. Furthermore, the common flux can be computed by adding contributions from the convective and diffusive components, which we denote with superscripts $(c)$ and $(v)$, respectively. Then, we can write
\begin{equation}
  \hat{{\bm{\mathfrak{F}}}}_{k,f} = \hat{{\bm{\mathfrak{F}}}}_{k,f}^{(c)} + \hat{{\bm{\mathfrak{F}}}}_{k,f}^{(v)},
\end{equation}
where we consider the following form of the common fluxes
\begin{align}
  \hat{{\bm{\mathfrak{F}}}}_{k,f}^{(c)} &= \bm F^{(c)}(\hat{u}^h_{\bar f}) + s^{(c)} (u^h_{k,f} - \hat{u}^h_{\bar f}) {{\bm n}}_{k,f},\\
  \hat{{\bm{\mathfrak{F}}}}_{k,f}^{(v)} &= \bm F^{(v)}(\hat{u}^h_{\bar f}, \bm q_{k, f}^{h}) + s^{(v)} (u^h_{k,f} - \hat{u}^h_{\bar f}) {{\bm n}}_{k,f},
  \label{eq:98hruon}
\end{align}
and $s^{(c)}$ and $s^{(v)}$ are convective and viscous stabilization parameters. Note that with this form of the fluxes, the definition of the common flux is given for each element, where the only information available is within itself and the trace variable. Conservation is implicitly enforced via transmission conditions
\begin{equation}
  \llbracket \hat{\bm{\mathfrak F}} \rrbracket_{\varepsilon^h_0} = 0,
    \label{eq:transmissionstrong}
\end{equation}
which can be discretely written 
\begin{equation}
  \sum_{\bar f\in \varepsilon^h_0} \int_{\bar f} \llbracket \hat{\bm{\mathfrak{F}}}\rrbracket_{\bar f} \phi ds+ \sum_{\bar f\in \varepsilon^h_\partial}  \int_{\bar f}  {\mathfrak{F}}^{\text{BC}}_{\bar f} \phi ds = 0,
  \label{eq:5treerf}
\end{equation}
and provides closure to the system. In these equations, ${\mathfrak{F}}^{\text{BC}}_{\bar f}$ is the normal boundary flux and the jump operator is defined at an interior face by 
\begin{equation}
  \llbracket\hat{\bm {\mathfrak{F}}} \rrbracket_{\bar f} = \hat{\bm {\mathfrak{F}}}_{k^+,f^+} \cdot { {\bm{n}}}_{k^+,f^+} + \hat{\bm {\mathfrak{F}}}_{k^-,f^-} \cdot {\bm{n}}_{k^-,f^-}.
  \label{eq:53904jn}
\end{equation}
After summing over all elements, we can state the hybridized form of the flux reconstruction approach for convection-diffusion type problems as follows
\begin{subequations}\begin{align}
  \sum_{\Omega_k\in \mathcal{T}_h}\tilde{\bm q}_k^h - \sum_{i=1}^{N_s} U_{k,i} \tilde\nabla \varphi_i(\tilde{\bm x}) - \sum_{f=1}^{N_f} \sum_{m=1}^{N_r} \tilde{\bm n}_f^m\cdot \tilde \nabla \cdot \bm g_f^m(\tilde{\bm x})\left[\mathfrak{U}_{k,f} - u_{k,f}^h\right]_{\tilde{\bm x} = \tilde{\bm x}^r_{f,m}} &= 0,\\
 \textcolor{\colornameedits}{ \sum_{\Omega_k\in \mathcal{T}_h}\frac{\partial {u}^h_k}{\partial t} + \frac{1}{J_k}\sum_{i=1}^{N_s} \tilde{\bm F}_{k,i} \cdot \tilde \nabla \varphi_i(\tilde{\bm x}) + \frac{1}{J_k}\sum_{f=1}^{N_f}\sum_{m=1}^{N_r} \tilde\nabla\cdot \bm g_f^m(\tilde{\bm x}) \left[{\tilde{H}}(\tilde{\bm x})_{k,f}\right]_{\tilde{\bm x} = \tilde{\bm x}^r_{f,m}}}&= 0,\label{eq:hfr1}\\
  \sum_{\bar f\in \varepsilon^h_0} \int_{\bar f} \llbracket \hat{\bm{\mathfrak{F}}}\rrbracket_{\bar f} \phi d\bar f+ \sum_{\bar f\in \varepsilon^h_\partial}   \int_{\bar f} {\mathfrak{F}}^{\text{BC}}_{\bar f} \phi d\bar f & = 0,\label{eq:hfr2}
\end{align}\label{eq:hfr}\end{subequations}
where we have readily taken the divergence of the flux and its correction to arrive at~\eqref{eq:hfr1}. Typically, hybridized methods make use of discontinuous or globally continuous function spaces for the trace variable, which can be respectively defined by
\begin{subequations}
\begin{align}
  M^h_p &= \{\mu \in L_2(\varepsilon^h)~:~\mu|_{\bar f} \in \mathbb{P}^p(\bar f),~\forall \bar f \in \varepsilon^h\},\\
  \bar M^h_p &= M^h_p \cap C^0(\varepsilon^h).
\end{align}\label{eq:hfrspaces}\end{subequations}
These finite-element spaces lead to the so-called HFR and a hybridized method with a smaller space known as EFR, as previously defined by Pereira and Vermeire in~\cite{pereiraPerformanceAccuracyHybridized2022}, in line with existing naming conventions of the HDG methods of Cockburn et al.~\cite{cockburnUnifiedHybridizationDiscontinuous2009}. \textcolor{\colornameedits}{Note that EFR schemes are a subset of HFR methods with a smaller trace basis function space. To distinguish methods with discontinuous traces from the global family of methods, we will refer to these as HFR throughout this work, and will use the terminology \textit{hybridized FR methods} when both continuous and discontinuous traces will be considered.}

\begin{figure}[htb]
  \centering
  \begin{subfigure}[b]{0.49\textwidth}
    \includegraphics[width=\textwidth]{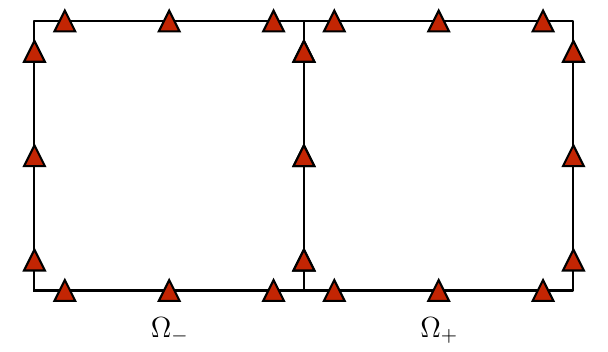}
    \caption{HFR}
  \end{subfigure}
  \begin{subfigure}[b]{0.49\textwidth}
    \includegraphics[width=\textwidth]{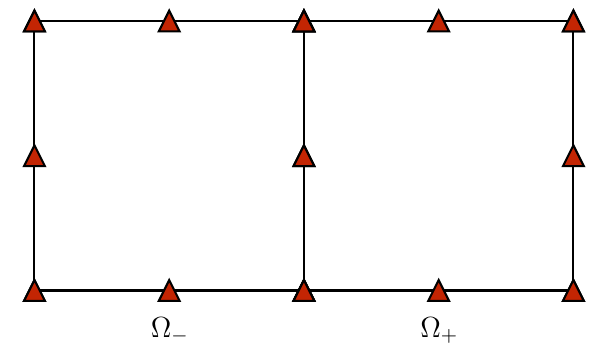}
    \caption{EFR}
  \end{subfigure}
  \caption[]{\textcolor{\colornameedits}{Trace variable location in an HFR (left) and EFR (right) discretization considering a $p=2$ scheme on the skeleton of two neighbouring quadrilateral elements}}
  \label{fig:tracepoints}
\end{figure}

\section{The Global System}

\subsection{Block Formulation}
Hybridized methods for linear advection-diffusion can be written in matrix form as follows
\begin{align}
\begin{bmatrix}
 \bm A & \bm B \\ 
 \bm C & \bm D 
\end{bmatrix}
\begin{bmatrix}
\bm u \\ 
\hat{\bm u}
\end{bmatrix}
= 
\begin{bmatrix}
\bm r \\ \bm s
\end{bmatrix},
\end{align}
where $\bm u$ and $\hat{\bm u}$ refer to the vector of internal and trace solution points, respectively. Because of the discontinuous nature of the interior solution and the decoupling resulting from the definition of the Riemann solvers, we can reduce the problem via static condensation and solve the condensed problem
\begin{equation}
  \bm L \hat{\bm u} = \bm t,
\end{equation}
where $\bm L = \bm D - \bm C \bm A^{-1} \bm B$ and $\bm t = \bm s - \bm C \bm A^{-1} \bm r$. Then, the solution can be obtained from
\begin{equation}
  \bm u = \bm A^{-1}(\bm r - \bm B \hat {\bm u}).
  \label{eq:localproblems}
\end{equation}
These matrices can be built efficiently per element
\begin{align}
    \label{eq:assembly} \bm L_{i,j} &= \bm L_{i,j} + \bm L^k_{\bar i,\bar j}, \\
    \label{eq:assembly2} \bm t_{i} &= \bm t_{i} + \bm t^k_{\bar i},
\end{align}
where the elemental matrices $\bm L^k$ and $\bm t^k$ are defined by
\begin{align}
   \bm L^k &:= \bm D^k - \bm C^k(\bm {A}^{-1})^k\bm B^k, \\
   \bm t^k &:= \bm s^k - \bm C^k (\bm A^{-1})^k\bm r^k,
\end{align}
and the indices $\bar i,~\bar j$ are associated with a surjective mapping of the element's flux points to the global trace points. The elemental blocks $\bm A^k$ and $\bm B^k$ are defined by
\begin{subequations}
\begin{align}
    \label{eq:elemA} &\bm A^k_{i,j} = \frac{1}{J_{k,i}}\left[\sum_{g=1}^{N_s} \tilde \nabla \varphi_g(\tilde{\bm x}) \cdot  \frac{\partial\tilde{\bm F}_{k,g}}{\partial u_{k,j}}  + \sum_{f=1}^{N_f}\sum_{m=1}^{N_r} \tilde\nabla\cdot \bm g_f^m (\tilde{\bm x}) \frac{\partial \tilde H_{k,f}}{\partial u_{k,j}}(\tilde{\bm x}^r_{f,m})\right]_{\tilde{\bm x} = \tilde{\bm x}_i^s},\\ 
    &\bm B^k_{i,(f-1)N_r+m} =  \frac{1}{J_{k,i}} \tilde\nabla\cdot \bm g_f^m (\tilde{\bm x}^s_i) \frac{\partial \tilde H_{k,f}}{\partial \hat u_{(f-1)N_r+m}}(\tilde{\bm x}^r_{f,m}),  
\end{align}\label{eq:elemsystem}\end{subequations}
where $i,j=1,\ldots,N_s$, $f=1\ldots,N_f$, $m=1,\ldots,N_r$. $u_{k,j}$ is the solution at the $j$-th solution point in element $\Omega_k$ and $\hat u_l$ is the trace value overlapping the $l$-th flux point. $\bm C^k$ and $\bm D^k$ can be trivially assembled from
\begin{subequations}
\begin{align}
  & \bar{\bm C}^k_{f,qj} = \frac{\partial \tilde{\hat{\bm{\mathfrak{F}}}}_{k,f}}{\partial u_{k,j}} (\tilde {\bm x}^r_{f,q}) \cdot \tilde{\bm n}_f^q, \\
  & \bar{\bm D}^k_{f,qt} = \frac{\partial \tilde{\hat{\bm{\mathfrak{F}}}}_{k,f}}{\partial \hat u_{k,f,t}} (\tilde {\bm x}^r_{f,q}) \cdot \tilde{\bm n}_f^q,
\end{align}\label{eq:elemsystemB}\end{subequations}
such that elementwise matrices
\begin{subequations}
\begin{align}
  &\bm C^k = (\bm M^k_0 \bar{\bm C}^k_0,\ldots,\bm M^k_{N_f} \bar{\bm C}^k_{N_f})^T,\\
  &\bm D^k = \operatorname{diag} (\bm M^k_0 \bar{\bm D}^k_0,\ldots,\bm M^k_{N_f} \bar{\bm D}^k_{N_f}),
\end{align}\end{subequations}
are obtained. In addition, $\bm M_k^{f}$ is the face mass matrix $\bm M^f_{k,qt}=\int_{\bar f}{\phi_q\phi_t}d\bar f$. Similarly, the vectors $\bm r^k,~\bm s^{k}$ evaluate the right-hand side of Equations~\eqref{eq:hfr} typically containing boundary data. In addition, the flux Jacobian blocks account for the sensitivity with respect to the gradient and the solution. Note that we have chosen to directly formulate the blocks in terms of $u,~\hat u$, but it is also possible to formulate in terms of $\bm q,u,\hat u$ and then perform the static condensation procedure. For a description of how to proceed for nonlinear problems and to include the effects of the temporal scheme, see~\cite{pereiraPerformanceAccuracyHybridized2022}.
\section{Stability Analysis}
The hybridized methods considered in~\cite{pereiraPerformanceAccuracyHybridized2022} were shown to recover conventional FR formulations for linear advection problems and that EFR methods introduced additional dissipation via spectral analysis. In this section, we analyze the behaviour of hybridized FR methods for linear-diffusion problems that make use of the Vincent-Castonguay-Jameson-Huynh (VCJH) correction functions~\cite{vincentNewClassHighOrder2011} and discuss connections to HDG methods for which linear stability proofs have been obtained. For this purpose, we consider the linear advection-diffusion equation
\begin{equation}
   \frac{\partial u}{\partial t} + \bm \alpha \cdot \nabla u - \beta \nabla^2 u = 0,
\end{equation}
where $\bm \alpha$ is the advection speed, and $\beta$ is the diffusion coefficient. We can rewrite this second-order problem as a system of first-order equations to be consistent with~\eqref{eq:systempde}
\begin{align}
   & \frac{\partial u}{\partial t} + \nabla \cdot (\bm \alpha u - \beta \bm q) = 0,\\
   &\bm q - \nabla u = 0,
\end{align}
subjected to periodic boundary conditions.
\subsection{Explicit forms of the numerical trace}\label{sec:explicittrace}
As previously stated in the implementation section, we consider hybridized forms of FR methods where the trace variable may belong to one of the finite-element spaces in~\eqref{eq:hfrspaces}, 
which lead to the so-called hybridized and embedded flux reconstruction schemes. At a given flux point, the Riemann flux for advection-diffusion is given by
\begin{equation}
  \hat{{\bm{\mathfrak{F}}}}_{k,f} = \bm \alpha \hat{u}^h_{\bar f} - \beta\bm q_{k,f}^h + s_{k,f} (u^h_{k,f} - \hat{u}^h_{\bar f}) {{\bm n}}_{k,f},
\label{eq:riemflux}
\end{equation}
where $s_{k,f}=s^{(c)}_{k,f}+s^{(v)}_{k,f}$. Because of the discontinuous nature of the trace polynomials in HFR with space $M^h$, the conservation condition is applicable pointwise. Application of the transmission conditions on an interior trace point to solve for $\hat u$ yields the following explicit expressions for the trace variable for the HFR method
\begin{align}
\hat u^{{M}_p^h} &= \frac{\ldcbrace s u\rdcbrace}{\ldcbrace s \rdcbrace} - \frac{\beta}{2} \frac{\llbracket{\bm q}\rrbracket}{\ldcbrace s \rdcbrace}.
\label{eq:explicittracehfr}
\end{align}
In the case of the EFR method, a simple explicit expression cannot be obtained solely in terms of $u$ due to the global coupling resulting from the reduced space of the trace. If, for instance, we choose to under-integrate the transmission conditions by employing a quadrature such as Gauss-Lobatto-Legendre (GLL), the interior solution coupling is reduced, and a simplified expression can be written as
\begin{equation}
  \hat u^{\bar{{M}}_p^{h,GLL}}  =  \frac{\sum_{\bar F}\ldcbrace s w u\rdcbrace_{\bar F}}{\sum_{\bar F}\ldcbrace s \rdcbrace_{\bar F}} - \frac{\beta}{2} \frac{\sum_{\bar F}\llbracket{ w \bm q}\rrbracket_{\bar F}}{\sum_{\bar F}\ldcbrace s \rdcbrace_{\bar F}}\label{eq:explicittraceefr},
\end{equation}
where $\bar F$ refers to the faces intersecting the trace point. Here $w$ is a quadrature weight arising from the diagonal local mass matrix of the GLL quadrature. However, in the rest of the work, we employ exact integration to mitigate possible aliasing issues of the GLL quadrature choice. The averaging and the jump operators are defined as follows
\begin{equation}
  \ldcbrace u \rdcbrace = \frac{u_- + u_+}{2}, \quad \quad \llbracket \bm q \rrbracket = \bm q_-\bm n_- + \bm q_+\bm n_+.
\end{equation}
From these expressions, assuming a homogeneous definition of the stabilization on the $(-)$ and $(+)$ sides of the interface for all faces, we note that the problem is undefined for $\ldcbrace s \rdcbrace=0$, and hence the following statement is a constraint for hybridized advection and advection-diffusion problems
\begin{equation}
  \ldcbrace s \rdcbrace\neq 0.
  \label{eq:sneq0}
\end{equation}
\subsection{Correction Functions}
The correction functions we consider in this work are those that define energy-stable FR (ESFR) methods in one dimension, i.e. the VCJH correction functions. They can be defined as follows
\begin{equation}
g_L = \frac{(-1)^p}{2} \left[L_p - \frac{\eta_p L_{p-1} + L_{p+1}}{1 + \eta_p}\right],\quad
g_R = \frac{1}{2} \left[L_p + \frac{\eta_p L_{p-1} + L_{p+1}}{1+\eta_p}\right],
\end{equation}
where $L_p$ is a Legendre polynomial of degree $p$ and 
\begin{equation}
  \eta_p = \frac{c(2p + 1)(a_p p!)^2}{2}, 
\end{equation}
with $a_p$ the leading coefficient in $L_p$, and $c$ a free parameter which can recover existing methods. For instance, $c_{DG}=c=0$ recovers the discontinuous Galerkin method, $c_{SD}$ recovers the spectral difference method. For completeness, we also consider the $c_{HU}$ method, as defined in~\cite{huynhFluxReconstructionApproach2007}.  These functions can be directly extended to tensor product elements. See~\cite{vincentNewClassHighOrder2011} for the numeric value of $c$ for these methods. 

\subsection{Proof of stability}
We now devise a methodology to show energy stability of hybridized flux reconstruction methods. This will provide insights into the stabilization mechanisms of hybridized methods alongside VCJH correction functions. Specifically, we consider the analysis of hybridized methods on quadrilateral Cartesian grids, which have transformation Jacobian matrices of the form
\begin{equation}
  \bm J_{k} = \begin{bmatrix} J_x & 0 \\ 0 & J_y  \end{bmatrix},
\end{equation}
and hence for these problems, the face Jacobians $\hat J_{\bar f}=J_{x}$ and $\hat J_{\bar f}=J_{y}$ at a horizontal and vertical face, respectively, since the cross-terms are zero. This analysis has been widely studied for conventional FR schemes with several advective and diffusive Riemann solvers for one~\cite{castonguayEnergyStableFlux2013,vincentNewClassHighOrder2011,quaegebeurStabilityEnergyStable2019} and higher-dimensional~\cite{quaegebeurStabilityEnergyStable2020,sheshadriErratumStabilityFlux2016} problems.

We seek to study the time evolution of the solution using a suitable Sobolev norm. We perform different algebraic manipulations to determine the constraints that will yield well-defined and linearly stable schemes. We make direct use of the proofs in the work of Sheshadri et al.~\cite{sheshadriErratumStabilityFlux2016,sheshadriAnalysisStabilityFlux2018,sheshadriStabilityFluxReconstruction2016} for conventional FR and augment them with algebraic manipulated forms of the transmission conditions. 

First, we introduce two important equations that will enable the study via the following lemmas.
\begin{lemma}\label{lemma:integration}
For hybridizable FR methods on Cartesian grids, the following holds
\begin{equation}
\sum_{k=1}^N \int\limits_{\partial\Omega_k} \hat{\bm{\mathfrak{F}}}_k\cdot\bm n_k \hat u^h ds= 0. 
\label{eq:addstabint}
\end{equation}
\end{lemma}
\begin{proof}
On multiplying the transmission conditions from Equation~\eqref{eq:transmissionstrong} by a function $\mu$ that belongs to a finite-dimensional space from~\eqref{eq:hfrspaces} and integrating over each $\bar f$, we have
\begin{equation}
  \sum_{\bar f\in \varepsilon^h} \int\limits_{\bar f} \left( \llbracket \hat{\bm{\mathfrak{F}}}\rrbracket_{\bar f} \right)\mu ds = 0,
\end{equation}
which can be rewritten using the identity 
\begin{equation}
  \sum_{\bar f\in \varepsilon^h} \int\limits_{\bar f} \left( \llbracket \hat{\bm{\mathfrak{F}}}\rrbracket_{\bar f} \right)\mu ds = \sum_{k=1}^N \int\limits_{\partial\Omega_k} \hat{\bm{\mathfrak{F}}}_k\cdot\bm n_k \mu ds = 0.
  \label{eq:identity}
\end{equation}
Recall that the problem has been defined to be periodic and that $\hat u^h$ also belongs to one of the spaces in~\eqref{eq:hfrspaces}. Since $\int_f \mu df= \int_f \mu|_fdf$ and both $\mu,~\hat u^h$ belong to the same space, we substitute $\mu$ by the trace variable $\hat u^h$ and the proof is complete. 
\end{proof}
\begin{lemma}\label{lemma:differentiation}
For all hybridizable FR methods, the conservativity condition holds and therefore, the following holds as well
\begin{equation}
\sum_{k=1}^N \sum\limits_f^{N_f}\left[J_{\psi_k}^{2p+1} \frac{\partial^p (\hat{\bm{\mathfrak{F}}}_k\cdot\bm n_k)}{\partial\psi^p}\frac{\partial^p \hat u^h}{\partial\psi^p}\right]_f = 0.
\label{eq:addstabdiff}
\end{equation}
\end{lemma}{}
\begin{proof}
The procedure is similar to the previous lemma but uses differentiation. On differentiating the transmission conditions $p$ times along the reference $\tilde \psi$-direction corresponding to each face, multiplying by $\frac{\partial^p \mu}{\partial\psi^p}$ and integrating over $\varepsilon^h$
\begin{equation}
  \sum_{\bar f\in \varepsilon^h} \int\limits_{\bar f} \frac{\partial^p ( \llbracket \hat{\bm{\mathfrak{F}}}\rrbracket_{\bar f})}{\partial\psi^p}\frac{\partial^p \mu}{\partial\psi^p} d\bar f = 0,
\end{equation}
since both $\hat{\bm{\mathfrak{F}}}\cdot \bm n$ and $\mu|_{\bar f} \in \mathbb{P}_p$, the integrand is a constant and hence
\begin{equation}
  \sum_{\bar f\in \varepsilon^h} \left[\frac{\partial^p ( \llbracket \hat{\bm{\mathfrak{F}}}\rrbracket_{\bar f})}{\partial\psi^p}\frac{\partial^p \mu}{\partial\psi^p}\right]_{\bar f} = 0,
  \label{eq:addstabdiff22}
\end{equation}
where we have omitted the integration limits since we are considering a Cartesian grid with constant nonzero face Jacobians. Substituting $\mu$ by the trace variable and applying identity~\eqref{eq:identity} completes the above lemma. Note that for $J_x=J_y$, multiplying~\eqref{eq:addstabdiff22} by any constant Jacobian will not modify the result. For consistency with the forms of the norm that will be derived, this factor is $J_{\psi_k}^{2p+1}$. Another way to obtain this factor is by taking derivatives in reference space following the steps above and then transforming the resulting equations to physical space. 
\end{proof}
\begin{lemma}\label{lemma:originalplushfr}
For the tensor-product FR formulation with VCJH correction functions, the following holds
\begin{equation}
 \frac{1}{2} \frac{d}{dt} \Vert u^h \Vert_{p,2}^2 = - \beta\Vert \bm q\Vert^2 + \Theta^{\mathrm{FR}} + \hat \Theta,
 \label{eq:frstab}
\end{equation}
where 
\begin{equation}
\begin{split}
\Theta^{\mathrm{FR}} &= -  \sum_{k=1}^N \int\limits_{\partial\Omega_k} u^h_k (\hat{\bm{\mathfrak{F}}}_k \cdot \bm n) ds  \\
&+ \sum_{k=1}^N \frac{1}{2}  \int\limits_{\partial\Omega_k}u^h_k (\bm F^{h,(c)}_k \cdot \bm n)  ds \\ 
&-  \sum_{k=1}^N \int\limits_{\partial\Omega_k} (\textcolor{\colornameedits}{\hat u^h} - u^h_k)(\bm F^{h,(v)}_k \cdot \bm n)  ds\\  
&- c \sum_{k=1}^N\sum_{f}^{N_f} \left[J_{\psi_k}^{2p+1}\frac{\partial^p u^h_k}{\partial\psi^p}\frac{\partial^p}{\partial \psi^p} (\hat{\bm{\mathfrak{F}}}_{k,f} \cdot \bm n)\right]_f \\
&+ c  \sum_{k=1}^N\sum_{f}^{N_f} \left[J_{\psi_k}^{2p+1}\frac{1}{2}\frac{\partial^p u^h_k}{\partial\psi^p}\frac{\partial^p}{\partial \psi^p} (\bm F^{h,(c)}_{k,f} \cdot \bm n)\right]_f\\
&- c \sum_{k=1}^N\sum_{f}^{N_f} \left[J_{\psi_k}^{2p+1} \frac{\partial^p (\textcolor{\colornameedits}{\hat u^h}-u^h_k)}{\partial\psi^p}\frac{\partial^p}{\partial \psi^p} (\bm F^{h,(v)}_{k,f} \cdot \bm n)\right]_f \\
\end{split}
\label{eq:thetafr}
\end{equation}
and
\begin{equation}
\hat \Theta = \sum_{k=1}^N \int\limits_{\partial \Omega_k} \textcolor{\colornameedits}{\hat u}^h (\hat{\bm{\mathfrak{F}}}_k\cdot \bm n)  ds + c\sum_{k=1}^N\sum_{f}^{N_f} \left[J_{\psi_k}^{2p+1} \frac{\partial^p \hat u^h}{\partial \psi^p} \frac{\partial^p (\hat{\bm{\mathfrak{F}}}_k\cdot \bm n)}{\partial \psi^p}\right]_f
\label{eq:thetahfr}
\end{equation}
\end{lemma}
\begin{proof}
Here we have directly introduced the expressions obtained from the proof of stability of the FR method by Sheshadri et al.~\cite{sheshadriErratumStabilityFlux2016,sheshadriAnalysisStabilityFlux2018}. The reader can refer to it for the proof of such a statement. After algebraic manipulations and changes in the notation for the sake of consistency, a general expression for the stability of the FR method on Cartesian quadrilateral elements is given by
\begin{equation}
  \frac{1}{2} \frac{d}{dt} \Vert u^h \Vert^2 = - \beta\Vert \bm q\Vert^2 + \Theta^{\mathrm{FR}},
\end{equation}
where $\Theta^{\text{FR}}$ reads as in Equation~\eqref{eq:thetafr}. Here, $\psi$ is a dummy coordinate variable such that $\psi=x$ for horizontal faces and $\psi=y$ for vertical faces. In addition, $\bm F^{(c)}$ and $\bm F^{(v)}$ refer to the advective and diffusive fluxes, and $\hat{\bm{\mathfrak{F}}}_{k}\cdot \bm n$ is the total normal Riemann flux involving both advection and diffusion. Typically, explicit forms of the numerical trace $\hat u^h$ are used to derive these stability proofs. However, since hybridized FR methods implicitly define it, we leave it as a variable for this analysis. This also allows different function spaces for the trace to be considered. In this sense, we augment the above expression with Lemmas~\ref{lemma:integration} and~\ref{lemma:differentiation}, which do not affect the results since they add to zero. Hence, multiplying Equation~\eqref{eq:addstabdiff} by $c$ and adding it to~\eqref{eq:addstabint} completes the proof.
\end{proof}
\noindent With these tools, we are ready to state our theorem on the stability of hybridized FR methods.
\begin{theorem}
Using a tensor-product formulation of the hybridized FR methods with VCJH correction functions, the two-dimensional linear advection-diffusion equation with periodic boundary conditions on Cartesian quadrilateral elements, it can be shown that if
\begin{itemize}
\item The stabilization parameters $s_\pm=s_\pm^{(c)}+s_\pm^{(v)}$ are chosen such that $s_{\pm} > \frac{\bm \alpha\cdot \bm n_\pm}{2}$ and
\item The correction parameter satisfies $c\geq 0$,
\end{itemize}
then the following expression holds
\begin{equation}
\frac{1}{2} \frac{d}{dt} \Vert u^h \Vert_{p,2}^2 \leq 0,
\end{equation}
for a broken Sobolev norm of the solution given by
\begin{equation}
  \Vert u^h \Vert_{p,2}^2 = \sum_{k=1}^N \int_{\Omega_k} \left[(u_k^h)^2 + \frac{c}{2}\left(\left(\frac{\partial^p u_k^h}{\partial \xi^p}\right)^2 + \left(\frac{\partial^p u_k^h}{\partial \eta^p}\right)^2\right) + \frac{c^2}{4}\left(\frac{\partial^{2p} u_k^h}{\partial\xi^p\partial\eta^p}\right)^2\right]d\Omega_k.
\end{equation}
\end{theorem}
\begin{proof}
To state this proof, we can rewrite the equations in Lemma~\ref{lemma:originalplushfr} as a summation over all faces in the computational domain. Note that this is valid since all integrations and derivatives in the previous equations are performed over the borders of the elements. Hence, we consider one of these faces with either horizontal or vertical direction with left and right elements $\Omega_-$ and $\Omega_+$ and with outward unit normal vectors $\bm n_-$ and $\bm n_+$, respectively. Let us now expand each of the terms in these equations at a given face. The first term in Equation~\eqref{eq:thetafr} can be written as follows
\begin{align}
  \nonumber -& \int_{\bar f}\left( u_-^h \left[\textcolor{\colornameedits}{\bm F(\hat u^h_{\bar f}, \bm q_-^h)} \cdot \bm n_- + s_-(u_-^h - \hat u^h_{\bar f})\right] + 
    u_+^h \left[\textcolor{\colornameedits}{\bm F(\hat u^h_{\bar f}, \bm q_+^h)} \cdot \bm n_+ + s_+(u_+^h - \hat u^h_{\bar f})\right]\right)d\bar f \\
   =& -\int_{\bar f}\left( u_-^h[\bm \alpha \cdot \bm n_- \hat u^h_{\bar f} - \beta \bm q_-^h \cdot \bm n_- + s_-(u_-^h - \hat u^h_{\bar f})]
    + u_+^h[\bm \alpha \cdot \bm n_+ \hat u_{\bar f}^h - \beta \bm q_+^h \cdot \bm n_+ + s_+(u_+^h - \hat u_{\bar f}^h)]\right)d\bar f,
\end{align}
where we have expanded the definitions of the Riemann solver according to Equation~\eqref{eq:riemflux} and considered a total stabilization parameter $s=s^{(c)}+s^{(v)}$. Similarly, the second term involving the convective flux can be expanded for this face
\begin{align}
  &\nonumber\frac{1}{2} \int_{\bar f} \left( u_-^h \bm F^{(c)}(\hat u^h_{\bar f}) \cdot \bm n_- + u_+^h \bm F^{(c)}(\hat u^h_{\bar f}) \cdot \bm n_+\right) d\bar f \\
  &= \frac{1}{2} \int_{\bar f} \left( \bm \alpha \cdot \bm n_-  (u_-^{h})^2 + \bm \alpha \cdot \bm n_+ (u_+^{h})^2\right)d \bar f,
\end{align}
and the third term involving the diffusion component of the flux can be written 
\begin{align}
  &- \nonumber \int_{\bar f} \left((\hat u^h_{\bar f} - u_-^h) \bm F^{(v)}(\bm q_-^h)\cdot \bm n_- + (\hat u^h_{\bar f}- u_+^h) \bm F^{(v)}(\bm q_+^h)\cdot \bm n_+\right)d\bar f \\
  &= - \int_{\bar f} \left((\hat u^h_{\bar f} - u_-^h) (- \beta \bm q_-^h\cdot \bm n_-) + (\hat u^h_{\bar f} - u_+^h) (- \beta \bm q_+^h \cdot \bm n_+)  \right) d\bar f.
\end{align}
Finally, we consider the first term of Equation~\eqref{eq:thetahfr}
\begin{align}
& \nonumber \int_{\bar f} \left(\hat u^h_{\bar f} [\bm \alpha \cdot \bm n_- \hat u^h_{\bar f} - \beta \bm q_-^h \cdot \bm n_- + s_-(u_-^h - \hat u^h_{\bar f})] + \hat u^h_{\bar f}[\bm \alpha \cdot \bm n_+ \hat u^h_{\bar f} - \beta \bm q_+^h \cdot \bm n_+ + s_+(u_+^h - \hat u^h_{\bar f})]\right)d\bar f\\
&= \int_{\bar f} \left(\hat u^h_{\bar f}[- \beta \bm q_-^h \cdot \bm n_- + s_-(u_-^h - \hat u^h_{\bar f})] + \hat u^h_{\bar f}[- \beta \bm q_+^h \cdot \bm n_+ + s_+(u_+^h - \hat u^h_{\bar f})]\right)d\bar f,
\end{align}
where we have used $\bm n_- = -\bm n_+$ to cancel out the advective flux on the trace variable. After adding all of the above contributions, we write
\begin{align}
  \Theta_{\bar f}^{A} = \int\limits_{\bar f} \left[\bar s_-(u_-^h - \hat u^h_{\bar f})^2 + \bar s_+(u_+^h - \hat u^h_{\bar f})^2\right] d\bar f,
  \label{eq:thetaa}
\end{align}
where we have introduced $\bar s_{\pm}=s_\pm - \frac{\bm \alpha \cdot \bm n_\pm}{2}$. Note the exchange of energy between the two adjacent elements is implicitly done via the trace variable.  In a similar manner, we can obtain the contributions from the derivative terms and write
\begin{equation}
\begin{split}
  \textcolor{\colornameedits}{\Theta_{\bar f}^{B} = J^{2p+1}_{\psi_k}\left[\bar s_-\left(\frac{\partial u_-^h}{\partial\psi^p} - \frac{\partial \hat u^h_{\bar f}}{\partial\psi^p}\right)^2 + \bar s_+\left(\frac{\partial u_+^h}{\partial\psi^p} - \frac{\partial \hat u^h_{\bar f}}{\partial\psi^p}\right)^2\right]_{\bar f}},
  \label{eq:thetab}
\end{split}
\end{equation}
for which we omit the derivation since it follows a similar procedure. Considering a periodic domain, the sum over all faces $\bar f\in\varepsilon^h$ results in the following stability statement for hybridized FR schemes on quadrilateral elements
\begin{equation}
\begin{split}
\frac{1}{2} \frac{d}{dt} \Vert u^h \Vert^2 
   =& - \beta \Vert \bm q\Vert^2  - \sum_{\bar f\in \varepsilon^h} \left(\Theta_{\bar f}^A + c\Theta_{\bar f}^B\right). 
\end{split}
\label{eq:stabstatement}
\end{equation}
From the above statements, we observe that for $c\geq 0$ and 
\begin{equation}
  \bar s_{\pm} \geq 0 \quad \Rightarrow \quad s_{\pm} \geq \frac{\bm \alpha \cdot \bm n_\pm}{2},
\end{equation}
the hybridized form of the FR for advection-diffusion satisfies
\begin{equation}
  \frac{d}{dt} \Vert u^h \Vert^2 \leq 0.
\end{equation}
However, from the explicit forms of the numerical trace defined in Equations~\eqref{eq:explicittracehfr} and~\eqref{eq:explicittraceefr}, we see that the method is undefined for $s_-=-s_+$ and hence the inequality becomes strict
\begin{equation}
  s_{\pm} > \frac{\bm \alpha \cdot \bm n_\pm}{2},
\end{equation}
or equally
\begin{equation}
  s_{k,f} > \frac{\bm \alpha \cdot \bm n_k}{2},
\end{equation}
where the strict inequality has to be satisfied on at least one face side~\cite{cockburnUnifiedHybridizationDiscontinuous2009}. Note that this proof of stability recovers that of the HDG method~\cite{nguyenImplicitHighorderHybridizable2009c} for $c=0$, and hence we show via this analysis that the stability proof of linear advection-diffusion HDG can be recovered from HFR methods with this particular value of $c$.
\end{proof}

\subsection{Connection to standard FR schemes}
In this section, we show the connection of hybridized methods with conventional FR formulations for a typical choice of the stabilization parameter. It is important to note that only discontinuous trace polynomials may recover existing FR formulations for purely convective problems. To establish a connection with standard FR formulations, we consider pure advection and pure diffusion scenarios. A typical choice of stabilization for problems involving advection and diffusion is 
\begin{equation}
  s_{k,f} =   s_{k,f}^{(c)} +  s_{k,f}^{(v)} = \lambda|\bm\alpha\cdot\bm n| + \tau,
  \label{eq:stabform}
\end{equation}
where $\lambda$ is an upwinding constant and $\tau$ is the so-called diffusion stabilization parameter. Several choices have been studied, but generally, one can take it to be
\begin{equation}
  \tau = \frac{\beta}{\ell},
\end{equation}
with $\ell$ a diffusive-length scale that results from dimensional analysis.
\subsubsection{Advection Regime}
First, we consider the case of pure advection ($\beta=0$). For the above choice of stabilization with $\tau=0$, the energy statement reads
\begin{equation}
\begin{split}
  \frac{1}{2} \frac{d}{dt} \Vert u^h \Vert^2 
   = 
   - \sum_{\bar f\in \varepsilon^h} & \left[\int_{\bar f} \frac{|\bm\alpha\cdot\bm n|}{2} 
   \left[\zeta_-(u_-^h - \hat u^h_{\bar f})^2 + \zeta_+(u_+^h - \hat u^h_{\bar f})^2 \right] d\bar f\vphantom{\left(\frac{\partial \hat u^h_{\bar f}}{\partial\psi^p}\right)^2}\right.\\
    +& \left.J_{\psi_k}^{2p+1}  \frac{|\bm\alpha\cdot\bm n|}{2}  c\left[\zeta_-\left(\frac{\partial u_-^h}{\partial\psi^p} - \frac{\partial \hat u^h_{\bar f}}{\partial\psi^p}\right)^2
  + \zeta_+\left(\frac{\partial u_+^h}{\partial\psi^p} - \frac{\partial \hat u^h_{\bar f}}{\partial\psi^p}\right)^2\right]_{\bar f}\right],
\end{split}
\label{eq:advproof}
\end{equation}
which implies that for the general case where the neighbouring interface solutions can take any arbitrary value, nonpositivity can be guaranteed for $\lambda\geq \frac{1}{2}$. Here we have defined 
\begin{equation}
  \zeta_\pm = 2\lambda\pm\operatorname{sign}(\bm\alpha\cdot\bm n_-).
\end{equation}
While this suggests that central-like approaches ($\lambda\rightarrow 0$) are not suitable choices for hybridization of pure advection problems, this can be mitigated by choice of discontinuous trace polynomials. As previously discussed in~\cite{pereiraPerformanceAccuracyHybridized2022,nguyenImplicitHighorderHybridizable2009c}, this finite-dimensional function space leads to the exact formulation of standard FR schemes, where the relationship with the trace variable and the interface solution value is equal when the same stabilization parameter is used on both sides, i.e.,
\begin{equation}
  (u_-^h - \hat u^h_{\bar f}) = (\hat u^h_{\bar f} - u_+^h),\quad s_-=s_+,
\end{equation}
which is a consequence of the pointwise validity of the transmission conditions for HFR. This results in a less strict range of stable $\lambda$ parameters, as the stability statement becomes
\begin{equation}
  \frac{1}{2} \frac{d}{dt} \Vert u^h \Vert^2 = \sum_{\bar f\in\varepsilon^{h}} \left[-\frac{\lambda}{2}\int\limits_{\bar f}|\bm\alpha\cdot\bm n|(u_- - u_+)^2 d\bar f
  - c\frac{\lambda}{2}\left[J_{\psi_k}^{2p+1} |\bm\alpha\cdot\bm n|\left(\frac{\partial u_-^h}{\partial\psi^p} - \frac{\partial u_+^h}{\partial\psi^p}\right)^2\right]_{\bar f}\right],
\end{equation}
consistent with the analysis of the FR methods in~\cite{sheshadriStabilityFluxReconstruction2016,sheshadriErratumStabilityFlux2016}. However, the implicit characteristic of $\hat u$ still requires $\lambda>0$ for the problem to be well-defined. This means that the exact central FR scheme for advection cannot be recovered with this type of stabilization. Note that we can find an explicit form of the linear-advection common flux with discontinuous trace polynomials and show that it takes the following forms~\cite{nguyenImplicitHighorderHybridizable2009c}
\begin{equation}
  \hat{{\bm{\mathfrak{F}}}}_\pm = \frac{\bm\alpha\cdot\bm n_- s_- + s_+ s_-}{s_+ + s_-} u_- + \frac{\bm\alpha\cdot\bm n_- s_+ - s_+ s_-}{s_+ + s_-} u_+.
  \label{eq:commonfluxexplicit}
\end{equation}
Hence, one possible way to define a central HFR method can be shown if the stabilization parameters are taken to be different on each side of the interface and take the following form
\begin{align}
  s_- &= \gamma|\bm\alpha\cdot\bm n| + \bm\alpha\cdot\bm n_-, \label{eq:hfrcentral1}\\
  s_+ &= f(\gamma,\bm \alpha) |\bm\alpha\cdot\bm n| - \bm\alpha\cdot \bm n_-,
  \quad f(\gamma,\bm\alpha) = \frac{\gamma\operatorname{sign}(\bm\alpha\cdot\bm n_-)}{2\gamma + \operatorname{sign}(\bm \alpha\cdot\bm n_-)},\label{eq:hfrcentral2}
\end{align}
for $\gamma\notin[\frac{1}{2},1],\gamma>0$. This recovers the exact central FR scheme for pure advection and no dissipative mechanism with a trace of the form
\begin{equation}
\hat u = \left(\frac{\operatorname{sign}(\bm \alpha\cdot\bm n_-)}{2\gamma} + 1\right)u_- - \frac{\operatorname{sign}(\bm\alpha\cdot\bm n_-)}{2\gamma} u_+.
\label{eq:hfrcentral3}
\end{equation}
The proof of~\eqref{eq:hfrcentral1}-\eqref{eq:hfrcentral3} can be shown by seeking the forms of the stabilization parameters in~\eqref{eq:commonfluxexplicit} that yield equal terms multiplying $u_-$ and $u_+$. However, it is well-known that fully-central methods are inconvenient for applications of physical interest due to their lack of dissipation.
\subsubsection{Diffusive Regime}
In the case of pure diffusion ($\bm \alpha=\bm0$), we set $\lambda=0$ in Equation~\eqref{eq:stabform} and obtain that the evolution of the $L_2$ energy is governed by
\begin{equation}
\begin{split}
  \frac{1}{2} \frac{d}{dt} \Vert u^h \Vert^2 
   = &- \beta \Vert \bm q\Vert^2 \\
   &- \tau \sum_{\bar f\in \varepsilon^h} \left[\int\limits_{\bar f} 
   ((u_-^h - \hat u^h_{\bar f})^2 + (u_+^h - \hat u^h_{\bar f})^2 ) d\bar f\right.\\
   & \quad+ J_{\psi_k}^{2p+1} c\left[\left(\frac{\partial u_-^h}{\partial\psi^p} - \frac{\partial \hat u^h_{\bar f}}{\partial\psi^p}\right)^2
  + \left(\frac{\partial u_+^h}{\partial\psi^p} - \frac{\partial \hat u^h_{\bar f}}{\partial\psi^p}\right)^2\right]_{\bar f}\left.\vphantom{\int\limits_{\bar f}} \right],
\end{split}
\label{eq:difproof}
\end{equation}
which shows that for this type of problem with an arbitrary positive diffusion coefficient $\beta$, stability is observed for any value of the viscous stabilization $\tau>0$. Contrary to the pure advection regime, hybridized methods for diffusion have a particular form of the numerical trace that cannot recover existing FR-LDG schemes for any finite value of $\tau$. LDG approaches make use of interface solution values that take the form
\begin{equation}
  \hat u^{\text{FR}} = \ldcbrace u \rdcbrace- \zeta \llbracket u \rrbracket,
  \label{eq:explicittraceldg}
\end{equation}
where $\zeta$ is a directional parameter. However, hybridized LDG (LDG-H) methods result in numerical traces defined by
\begin{equation}
  \hat u^{\text{HFR}} = \ldcbrace u \rdcbrace- \zeta \llbracket u\rrbracket - \beta\theta \llbracket \bm q \rrbracket.
  \label{eq:explicittraceldgh}
\end{equation}
Comparing the explicit definitions in Section~\ref{sec:explicittrace}, we see that 
\begin{equation}
  \theta = \frac{1}{2\ldcbrace s \rdcbrace},
\end{equation}
and hence $\theta\neq 0$ for any finite value of the stabilization. We also note that~\eqref{eq:difproof} is consistent with the HDG method in~\cite{nguyenImplicitHighorderHybridizable2009c} if we set $c=0$.  
\section{Local Post-processing}\label{sec:post-processing}
Post-processing techniques have been widely used to improve the accuracy of  numerical solutions. This is possible due to the optimal convergence rates of the solution and flux in locally-conservative methods involving diffusion operators~\cite{cockburnUnifiedHybridizationDiscontinuous2009,cockburnSuperconvergentDiscontinuousGalerkin2009a,cockburnSuperconvergentLDGhybridizableGalerkin2008}. Thus, this approach will only be applied to our problems with discontinuous trace variables and not the EFR method, whose DG equivalent has been shown to display suboptimal flux convergence in~\cite{cockburnAnalysisEmbeddedDiscontinuous2009} since the flux is not single-valued at the flux points. This approach has been applied to numerous types of problems involving steady-state and time-dependent problems. The procedure generally consists of using a Raviart-Thomas projection of the flux  to obtain a better approximation in $H(\text{div};\Omega)$ and solving a local problem for the solution. Recall that in the FR approach, we use correction functions to create a $C_0$-continuous flux function. With these correction functions, the discontinuous flux can be reconstructed to take the values of the Riemann fluxes, upgrading it to $\mathbb{P}_{p+1}$. In this section, we present a modified version of the post-processing method presented in~\cite{nguyenImplicitHighorderHybridizable2009c}, which leverages the operators already defined in the FR framework without the need to create additional RT formulations. 

The first step in obtaining an elementwise superconvergent solution $u^{h*}_k$ is to post-process the flux. While this first step can be done on the total advective and diffusive fluxes for linear problems, we choose to work with the viscous flux only for the sake of simplicity. On each element, we reconstruct the diffusive flux by computing
\begin{equation}
  \tilde {\bm F}_k^{*(v)} = \bm F_k^{hD(v)}({\tilde{\bm x})|_{\tilde{\bm x}_s^{*}}} + \sum_{f=1}^{N_f} \sum_{m=1}^{N_r^*} \bm g_f^{m}(\tilde{\bm x})[\tilde{ H}(\tilde{ \bm x})]_{\tilde{\bm x} = \tilde{\bm x}_f^{*m}},
\end{equation}
where $\tilde{\bm x}_f^{*m}$ are $N_s^*$ post-processing points that define polynomials of degree $p^*=p+1$ and $\bm g$ is the correction vector function of the same degree. \textcolor{\colornameedits}{$N_r^*$ is the number of face flux points in the new space}. From this post-processed flux, we can now obtain a more accurate solution by solving 
\begin{subequations}
\begin{align}
\begin{split}
&-\beta \left[\sum_{i=1}^{N_s^*} \tilde \nabla\varphi_i^* \cdot {\nabla} {u}^*_{k,i} -\!\! \sum_{f=1}^{N_f}\sum_{m=1}^{N_{r}^*} \tilde\nabla \cdot \bm g_f^{*m} ({\nabla} {u}_{k}^* \cdot \tilde{\bm n}_{f}^m) \Big|_{\bm x_{f,m}^r} \right] \!\!\\&= \!\!\sum_{i=1}^{N_s^*} \tilde{\bm F}_{k,i}^{*(v)} \cdot \tilde\nabla \varphi_i^* - \sum_{f=1}^{N_f}\sum_{m=1}^{N_{r}^*} \tilde\nabla \cdot \bm g_f^{*m} (\tilde {\bm F}_{k,f}^{*(v)} \cdot \tilde{\bm  n}_{f}^m)\Big|_{\bm x_{f,m}^r}, \end{split} \\
&\int_{\Omega_k}(u^h_k - u^{h*}_k)d\bm x = 0.
\end{align}
\end{subequations}
where 
\begin{equation}
  \nabla u^*_k = \bm J_k^{-T} \tilde \nabla u_k^*(\tilde {\bm x}) = \bm J_k^{-T}  \sum_{j=1}^{N_s^*}  U_{k,j}^*\tilde \nabla \varphi_j^*,
\end{equation}
which can be shown to be the FR discretization of the following problem at the element level
\begin{subequations}
\begin{align}
  \nabla \cdot (-\beta \nabla u) &= \nabla \cdot  {\bm F}_k^{*(v)}, \\
  -\beta \nabla u \cdot \bm n &= {\bm F}_k^{*(v)} \cdot \bm n,\\
  \int_{\Omega_k}(u - u^*)d\bm x &= 0.
\end{align}
\end{subequations}
The last statement ensures elementwise conservation of the solution. This post-processing leverages the existing FR correction functions to enable superior convergence of the methods and extends the post-processing schemes to the full family of FR schemes. This post-processing can be applied to schemes involving any of the VCJH correction functions and recovers a scaled form of the conventional HDG post-processing when $\bm g$ is constructed with $c_{DG}$. Otherwise, the post-processing schemes seem to be new. Later in this work, we perform numerical examples to showcase the superconvergent characteristics of hybridized FR methods. In all cases, we use the same correction function for both the solution and the post-processing steps for the sake of consistency.

\section{Numerical Examples}
In this section, we perform numerical experiments to discuss the stability, performance and accuracy of hybridized FR methods in advection-diffusion problems. We mainly consider three values of the $c$ parameter that recover existing high-order formulations. These include $c_{DG}$, $c_{SD}$ and $c_{HU}$ for $p=1$ to $p=4$ schemes. For the linear problems, the post-processing scheme is applied with the same values of these correction parameters. The $L_2$-norm of the solution error is measured via
\begin{equation}
  E_{L_2}(\Omega_h) = \sqrt{\frac{1}{|\Omega|}\sum_{k=1}^N \int_{\Omega_k} (u^h_k - u^e)^2 d\Omega_k},
\end{equation}
for verification purposes, where $u^{e}$ is the analytical solution. Furthermore, we first consider a linear steady-state problem and two unsteady cases. Then, we present a nonlinear problem involving the compressible Navier-Stokes equations. All simulations are carried out serially on a 3.2 GHz Intel Core i5-5600 processor with 16Gb of RAM. The implicit system makes use of an exact Jacobian and is solved via the default Block-Jacobi preconditioner in the PETSc framework~\cite{petsc-user-ref}. \textcolor{\colornameedits}{Unless otherwise noted, all timing quantities such as wall-clock times are indicated in seconds.} 
\subsection{Steady-State Linear Advection-Diffusion}
Consider the linear advection-diffusion equation with a source term chosen such that the exact solution is given by
\begin{equation}
  u(\bm x) = xy \frac{(1-e^{(x-1)a_x})(1-e^{(y-1)a_y})}{(1-e^{(1-a_x)})(1-e^{(1-a_y)})},
\end{equation}
defined on $\Omega=[0,1]^2$. Dirichlet boundary conditions are obtained directly from the exact solution. This case has been used to analyze the accuracy and post-processing of steady-state HDG methods for the weakly convection-dominated regime \textcolor{\colornameedits}{in~\cite{nguyenImplicitHighorderHybridizable2009c}}. The advection velocity is set to $\bm \alpha = [25, 25]$ and the diffusion coefficient to $\beta=1$. Due to the relatively large advection velocity, a boundary layer is expected to form toward the right and top ends of the domain. Hence, we consider the $L_2$-norm of the error in a reduced space $\Omega^{L_2}=[0.1,0.9]^2$ to exclude the resolution of the boundary layer. The grid was generated using the following stretching function
\begin{equation}
  \bm x = \frac{1}{a}\tanh\left(\frac{[i, j]}{\sqrt{N}}\arctanh a\right),\quad 0 \leq i,j \leq \sqrt{N}-1,
  \label{eq:stretchingfunc}
\end{equation}
where $N$ is the total number of elements and we set $a=0.95$.
\begin{figure}[htb]
\sbox0{\hwplotDGmarker}\sbox1{\hwplotSDmarker}\sbox2{\hwplotHUmarker}%
  \centering
  \begin{subfigure}[b]{0.35\textwidth}
    \includegraphics[width=\textwidth,trim={12cm 4cm 12cm 4cm},clip]{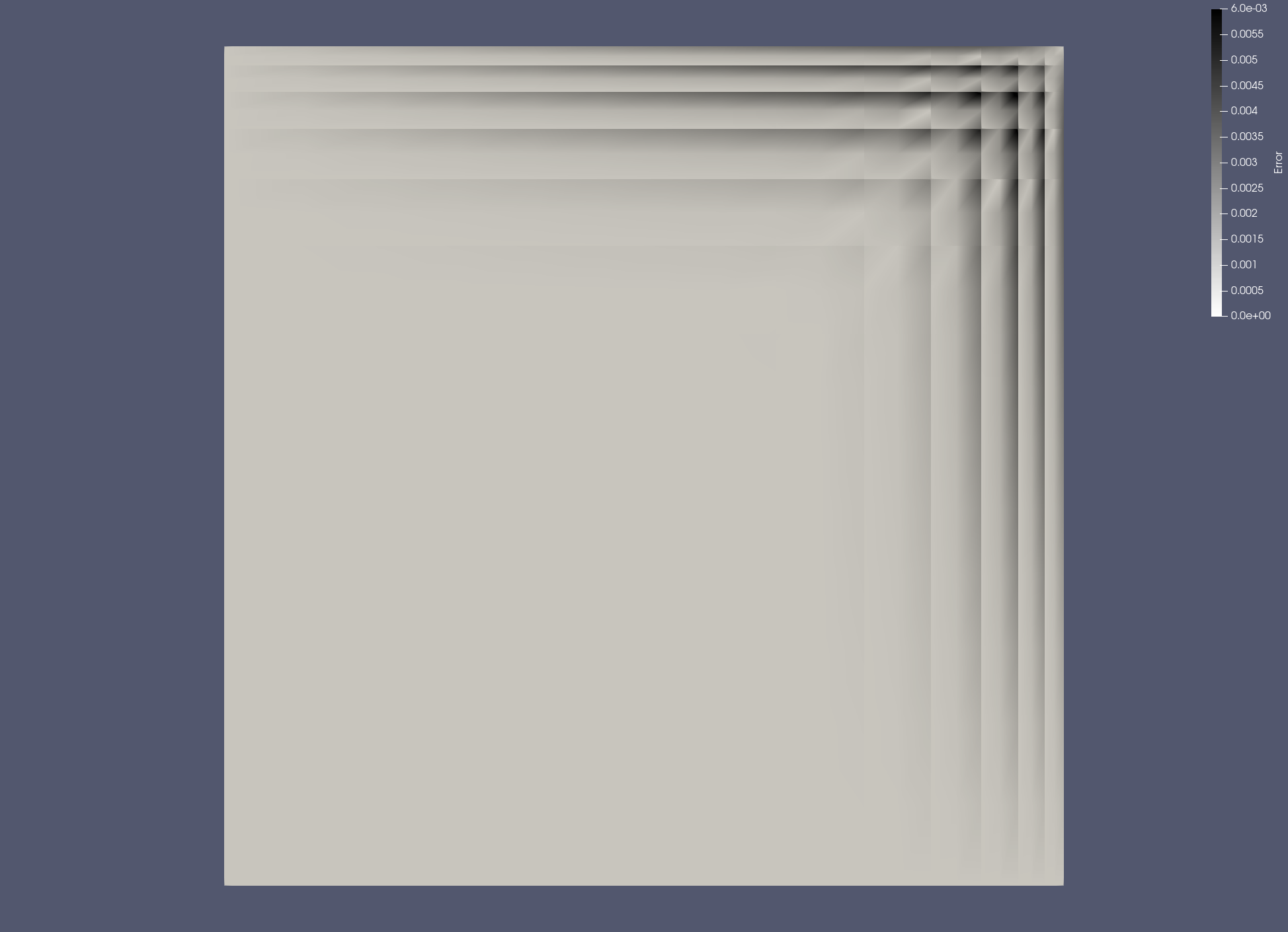}
    \caption{Error without post-processing}
  \end{subfigure}
    \begin{subfigure}[b]{0.35\textwidth}
    \includegraphics[width=\textwidth,trim={12cm 4cm 12cm 4cm},clip]{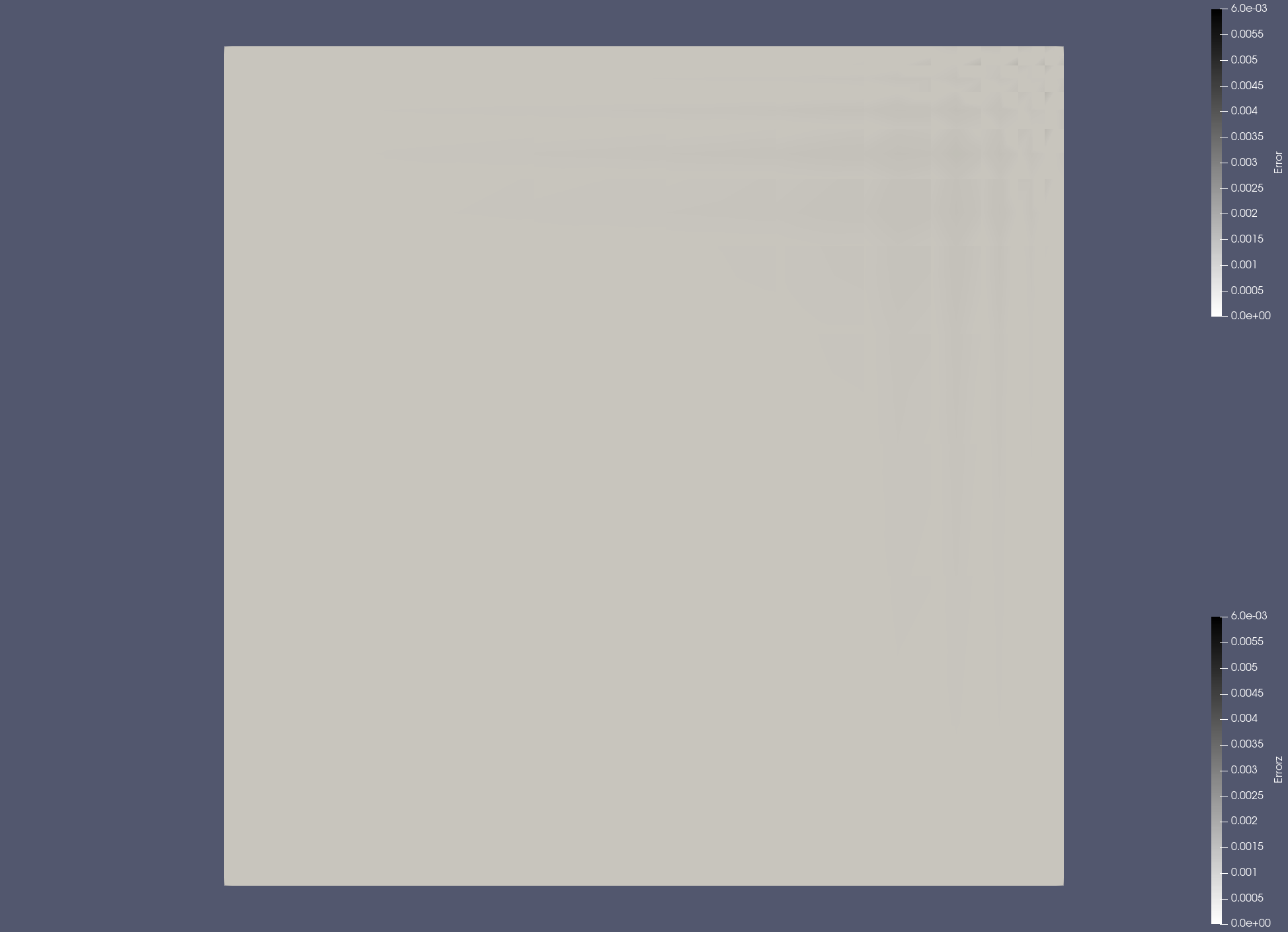}
    \caption{Error after post-processing}
  \end{subfigure}
  \caption{Contours of the solution and post-processed solution for the steady-state linear advection-diffusion problem for a $p=2$ solution on a $10\times10$ grid and $c=c_{SD}$. Linear color scaling adjusted to $[0,6\times 10^{-3}]$ from white to black}
  \label{fig:sadf-contour}
\end{figure}
\begin{figure}[htb]
\sbox0{\hwplotDGmarker}\sbox1{\hwplotSDmarker}\sbox2{\hwplotHUmarker}%
  \centering
  \begin{subfigure}[b]{0.49\textwidth}
    \includegraphics[width=\textwidth]{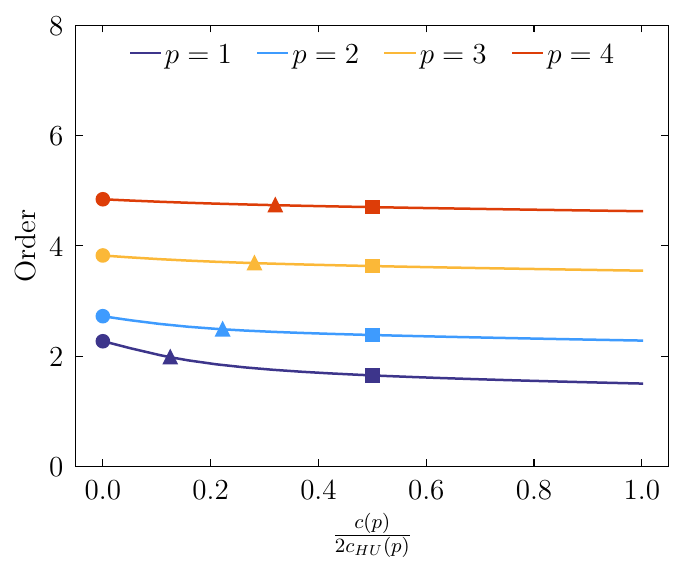}
    \caption{Solution}
  \end{subfigure}
    \begin{subfigure}[b]{0.49\textwidth}
    \includegraphics[width=\textwidth]{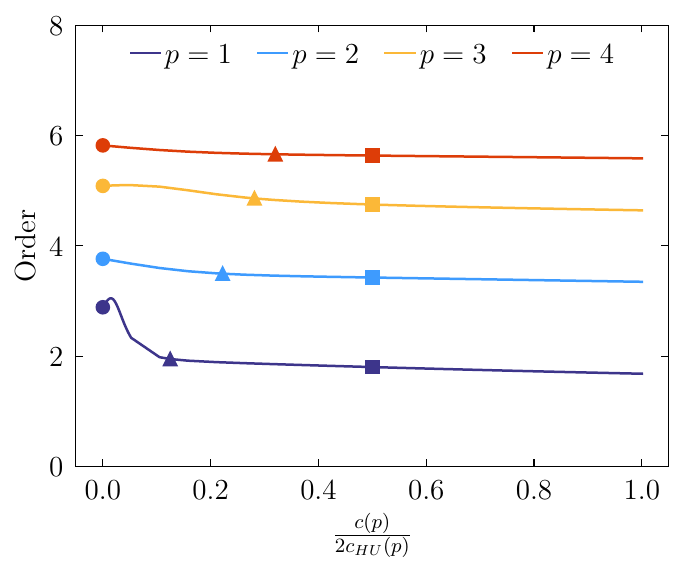}
    \caption{Post-processed solution}
  \end{subfigure}
  \caption{Effects of the correction function on the grid convergence rate of hybridized FR methods with discontinuous traces for the steady-state linear advection-diffusion problem. Markers for $c_{DG}$ (\usebox0), $c_{SD}$ (\usebox1), $c_{HU}$ (\usebox2) highlight these corrections at polynomial degrees $p=1$ to $p=4$}
  \label{fig:sadf-c-curved}
\end{figure}

We make use of four levels of refinement with $5\times5$, $10\times10$, $20\times20$ and $40\times40$ quadrilateral elements. Table~\ref{tab:ordersadvdiff} shows the $L_2$-norm of the error for $p=1$ to $p=4$ standard and hybridized FR schemes with correction parameters $c_{DG},~c_{SD}$ and $c_{HU}$. From the stability section, we observed that $c$ acts as an added dissipation mechanism to the $c_{DG}$ schemes, for which $c=0$. In this table, it can be observed that all schemes achieved the expected $p+1$ order of accuracy in logarithmic scale and that $c_{DG}$ has the smallest $L_2$ error for each of the considered types of discretization. Among these, the considered hybridized formulations in this work are more accurate than the standard FR-LDG formulation at the coarsest to finest levels. Specifically, the EFR method displayed the smallest error levels compared to HFR and FR for the finer levels of refinement. For instance, $p=3$ schemes with $c_{DG}$ displayed errors of $6.16\times10^{-8}$, $5.99\times10^{-8}$ and $5.01\times10^{-8}$ for FR, HFR and EFR, respectively. For HFR, where discontinuous traces are used, and for conventional FR methods, we applied the post-processing scheme in Section~\ref{sec:post-processing} and show the results obtained in Table~\ref{tab:orderzsadf}. As established in the literature, LDG methods are defined as in Equation~\eqref{eq:explicittraceldg}, and they do not possess the superconvergence property in the general case~\cite{cockburnSuperconvergentDiscontinuousGalerkin2009a}, but can be shown to superconverge for Cartesian grids with special choices of the common fluxes~\cite{cockburnSuperconvergenceLocalDiscontinuous2001}. For the sake of completeness, we show the results for this specific configuration in Table~\ref{tab:orderzsadvdiffconsistent} of this work but do not consider it anymore as it is not a feature of the arbitrary case. Interestingly, all $c_{DG}$ methods for HFR were able to achieve the expected $p+2$-order of accuracy after post-processing. However, when $c\neq c_{DG}$, only methods with $p>1$ were able to achieve the superconvergent behaviour. An example contour with the error in the computational domain is displayed in Figure~\ref{fig:sadf-contour}, where the error levels can be seen to decrease by an order of magnitude.
\begin{figure}[htb]
  \sbox0{\hwplotDGmarker}\sbox1{\hwplotSDmarker}\sbox2{\hwplotHUmarker}%
  \centering
  \begin{subfigure}[b]{0.49\textwidth}
    \includegraphics[width=\textwidth]{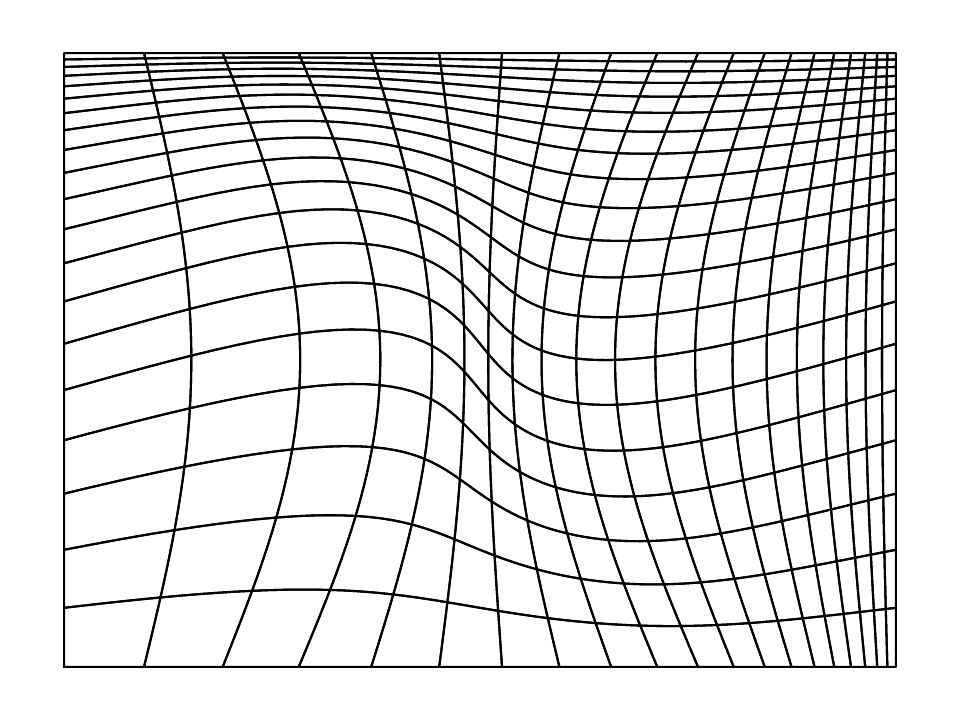}
  \end{subfigure}
  \caption{\textcolor{\colornameedits}{$20\times20$ curved element using $\bm x =\bm x_0 + A\sin(2\pi x)\sin(\pi y)$ with $A=0.1$ for the steady-state linear-advection diffusion problem. Here, $\bm x_0$ is obtained with Equation~\eqref{eq:stretchingfunc}}}
  \label{fig:uadf-curvedmesh}
\end{figure}

\begin{figure}[htb]
\sbox0{\hwplotDGmarker}\sbox1{\hwplotSDmarker}\sbox2{\hwplotHUmarker}%
  \centering
  \begin{subfigure}[b]{0.49\textwidth}
    \includegraphics[width=\textwidth]{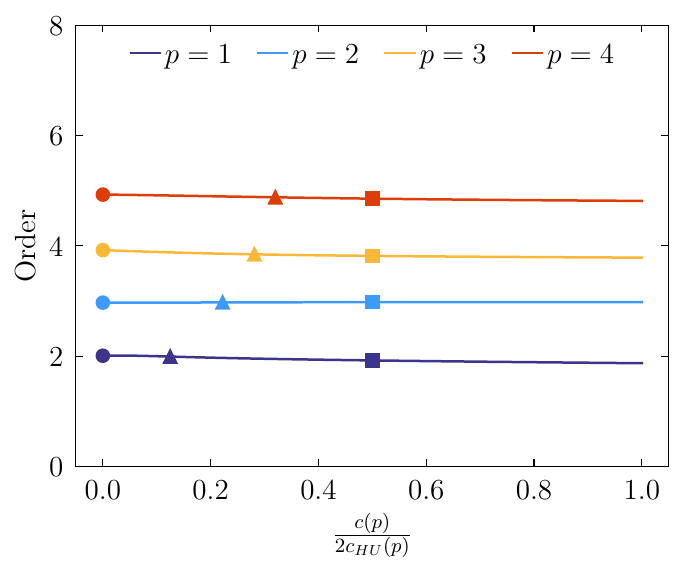}
    \caption{Solution}
  \end{subfigure}
    \begin{subfigure}[b]{0.49\textwidth}
    \includegraphics[width=\textwidth]{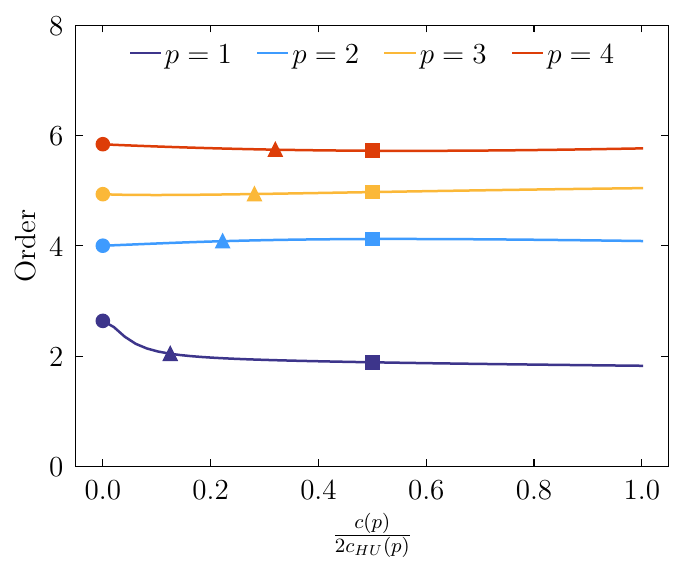}
    \caption{Post-processed solution}
  \end{subfigure}
  \caption{\textcolor{\colornameedits}{Effects of the correction function on the grid convergence rate of hybridized FR methods with discontinuous traces for the steady-state linear advection-diffusion problem on curved grids. Markers for $c_{DG}$ (\usebox0), $c_{SD}$ (\usebox1), $c_{HU}$ (\usebox2) highlight these corrections at polynomial degrees $p=1$ to $p=4$}}
  \label{fig:sadf-c}
\end{figure}
\textcolor{\colornameedits}{To further visualize the impact of the $c$-parameter on this superconvergent behaviour, we perform an additional set of simulations and compute the order of accuracy for values of $c\in[0,2c_{HU}]$ and display the results in Figure~\ref{fig:sadf-c}. It is known that the expected $p+1$ convergence of FR methods is lost at large values of $c$~\cite{vincentInsightsNeumannAnalysis2011}. For the considered relatively small range of $c$-parameters, the order of accuracy of the solution before post-processing slowly reduces as $c$ increases, but consistent with the results in the convergence study, rapid increase of the order is seen at $p=1$ from second to third order in the vicinity of $c\rightarrow0$. Similar results are obtained for the same problem using a curved grid with mapping degree 5, as shown in Figure~\ref{fig:uadf-curvedmesh}. The observed order of accuracy is consistent with $p+1$ and is obtained for all polynomial degrees. As the value of $c$ increases in these simulations, the convergence slowly reduces, but remains close to the expected value. Applying the post-processing scheme to the curved configurations showcases one additional order of convergence for $p>1$ for all values of $c$. Similarly, for $p=1$, the expected $p+2$ result is only observed in the vicinity of the $c_{DG}$ correction.}

\begin{figure}[htb]
\sbox0{\hwplotDGmarker}\sbox1{\hwplotSDmarker}\sbox2{\hwplotHUmarker}%
  \centering
  \begin{subfigure}[b]{0.49\textwidth}
    \includegraphics[width=\textwidth]{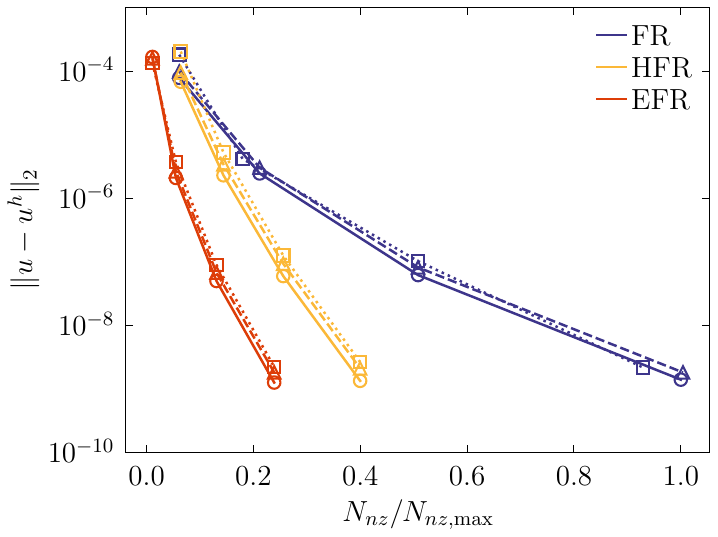}
    \caption{Relative number of nonzeros}
    \label{fig:sadf-metricsa}
  \end{subfigure}
    \begin{subfigure}[b]{0.49\textwidth}
    \includegraphics[width=\textwidth]{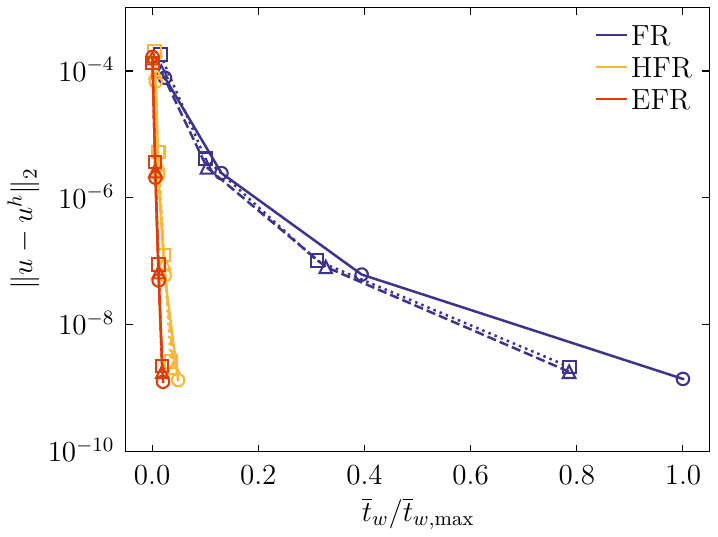}
    \caption{Relative solver time}
    \label{fig:sadf-metricsb}
  \end{subfigure}
  \begin{subfigure}[b]{0.49\textwidth}
    \includegraphics[width=\textwidth]{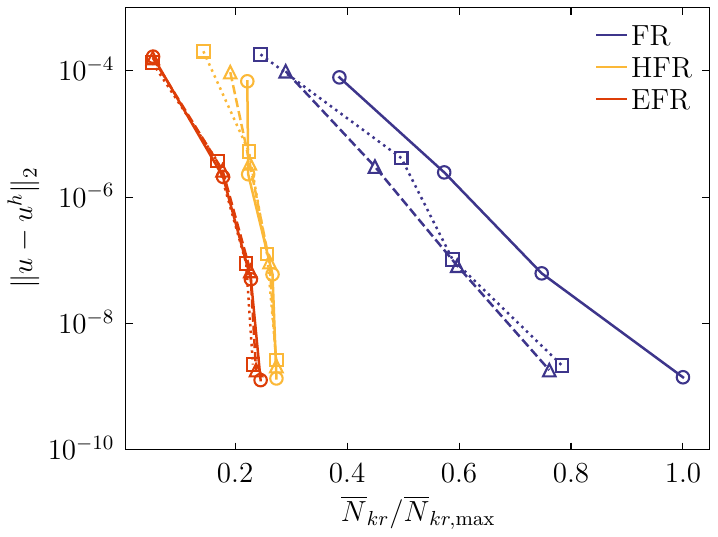}
    \caption{Number of implicit iterations}
    \label{fig:sadf-metricsc}
  \end{subfigure}
  \caption{Performance metrics for the steady-state linear advection-diffusion problem. Markers for $c_{DG}$ (\usebox0), $c_{SD}$ (\usebox1), $c_{HU}$ (\usebox2) have been added at polynomial degrees $p=1$ to $p=4$. For reference $N_{nz,\max}=1336250$, $\bar t_{w,max}=2.00s$, $\bar N_{kr,\max}=420$}
  \label{fig:sadf-metrics}
\end{figure}
Furthermore, we analyze the performance of using hybridized methods as opposed to implicit FR schemes in terms of the number of nonzeros, time spent per time step size and the number of GMRES iterations. We denote them $N_{nz}$, ${t}_w$ and $N_{kr}$, respectively. A bar is also placed on top of these quantities when we have considered their average per linear solve. For hybridized methods, $t_w$ accounts for both the solution of the system and the recovery of the internal solution via the local problems in Equation~\eqref{eq:localproblems}. Results are shown in Figure~\ref{fig:sadf-metrics} for each of the considered schemes and polynomial degrees on the finest $40\times40$ grid to reduce timing errors. It can be seen in Figure~\ref{fig:sadf-metricsa} that, for the highest polynomial degree considered $(p=4)$, the number of nonzeros in the system for HFR is reduced by half, and by about four times for EFR in comparison with FR schemes. Despite these metrics, it can be seen that the time spent on solving these systems before final convergence was achieved with a small fraction of that used in FR. Specifically, considering the $p=4$ simulations, we observe a reduction of 20-22 times for the HFR method and 40-48 times for the EFR method, depending on the value of $c$. This large difference can be attributed to the number of GMRES iterations shown in Figure~\ref{fig:sadf-metricsc}. FR schemes required a significantly larger number of implicit iterations to reach convergence for the Krylov solver, contrary to the hybridized schemes. Hence, hybridized FR methods have a significant benefit over standard implicit FR for all values of $c$ in terms of performance and accuracy for this steady-state problem.

\subsection{Advection-Diffusion of a Sine Wave}
Next, verification of unsteady linear advection-diffusion is performed in this section. This case was adapted from~\cite{castonguayEnergyStableFlux2013} to the two-dimensional setting. Consider the unsteady linear diffusion equation with initial condition
\begin{equation}
  u (x,y) = \sin \left(\frac{2\pi}{L} x\right)\sin\left(\frac{2\pi}{L} y\right),
\end{equation}
in a square of domain $L$ with periodic boundary conditions. The viscous stabilization was chosen to be $\tau=\beta=0.1$ and the simulation was run for one cycle on grids of $5\times5,~10\times10,~20\times20$ and $40\times40$. Then, the $L_2$ norm of the error was computed considering the analytical solution 
\begin{equation}  
  u^e(x,y, t) = \exp(-2\beta\pi^2t)\sin{[\pi(x-\alpha_xt)]}\sin{[\pi(y-\alpha_yt)]},
\end{equation}
after one cycle. The time discretization is done with a third-order single-diagonally implicit Runge Kutta method (SDIRK$_{3}$) using a \textcolor{\colornameedits}{consistent time step size $\Delta t^* =\Delta t/t_c = 5\times10^{-5}$ in all grids to reduce the temporal error. Here, $t_c$ is the convective characteristic time, which for this problem $t_c=\alpha_x/L=\alpha_y/L$ since we have chosen $\alpha_x=\alpha_y=1$}. Table~\ref{tab:orderuadf} shows the orders of accuracy for the considered correction functions and polynomial degrees 1 to 4, where $p+1$ convergence was observed for all considered schemes. Similar to the previous problem, hybridized methods displayed smaller error levels than the standard FR discretization for all considered values of $c$, with EFR being the most accurate at the finer grid levels. Specifically, $p=3$ schemes show errors in the $L_2$ norm of $2.07\times10^{-8}$, $1.25\times10^{-8}$ and $1.21\times10^{-8}$ for the FR, HFR and FR methods with $c_{DG}$, respectively. Hence, a reduction by about half of the error was obtained with the hybridized schemes. Note that, as before, $c_{DG}$ showed the best accuracy out of the three correction functions. 

The solution was post-processed using correction functions of the same $c$-parameter. Results and convergence orders are shown in Table~\ref{tab:orderzuadf} for the post-processed solution. FR-LDG did not show super accuracy for the case in which the switch was taken to be in an arbitrary fashion, but HFR achieved the expected order $p+2$ for the post-processed solution. Similar to the previous case, methods with $c\neq c_{DG}$ do not exhibit the $p+2$ order of accuracy at $p=1$.

\begin{figure}[htb]
  \sbox0{\hwplotDGmarker}\sbox1{\hwplotSDmarker}\sbox2{\hwplotHUmarker}%
  \centering
  \begin{subfigure}[b]{0.49\textwidth}
    \includegraphics[width=\textwidth]{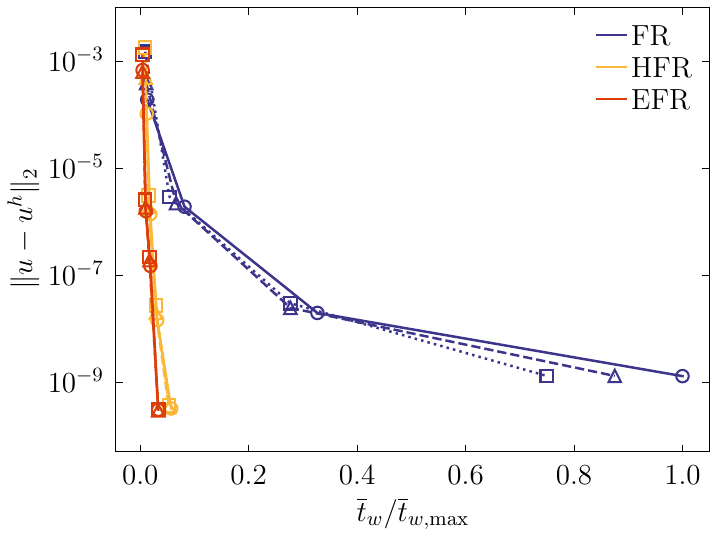}
    \caption{Relative solver time with $\Delta t^* = 2.5\times 10^{-4}$}
      \label{fig:uadf-metricsa}
  \end{subfigure}
    \begin{subfigure}[b]{0.49\textwidth}
    \includegraphics[width=\textwidth]{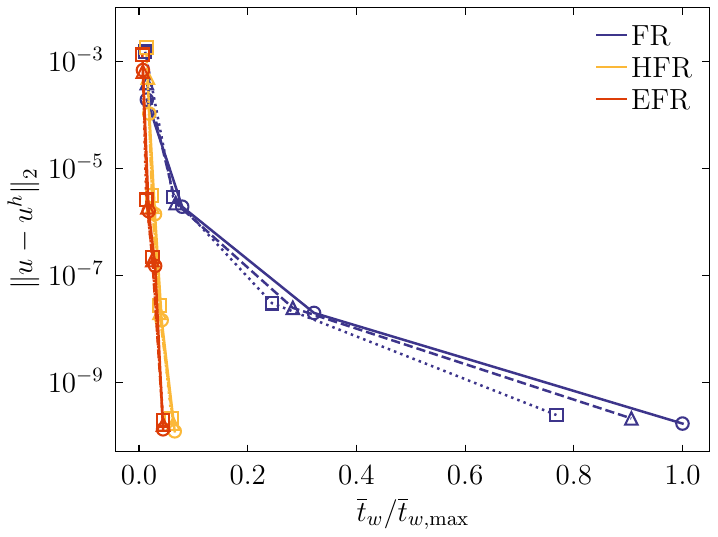}
    \caption{Relative solver time with $\Delta t^* = 5\times 10^{-5}$}
      \label{fig:uadf-metricsb}
    \end{subfigure}
  \caption{Performance metrics for the unsteady linear advection-diffusion problem. Markers for $c_{DG}$ (\usebox0), $c_{SD}$ (\usebox1), $c_{HU}$ (\usebox2) have been added at polynomial degrees $p=1$ to $p=4$.}
  \label{fig:uadf-metrics}
\end{figure}

\begin{figure}[htb]
  \sbox0{\hwplotDGmarker}\sbox1{\hwplotSDmarker}\sbox2{\hwplotHUmarker}%
  \centering
  \begin{subfigure}[b]{0.49\textwidth}
    \includegraphics[width=\textwidth]{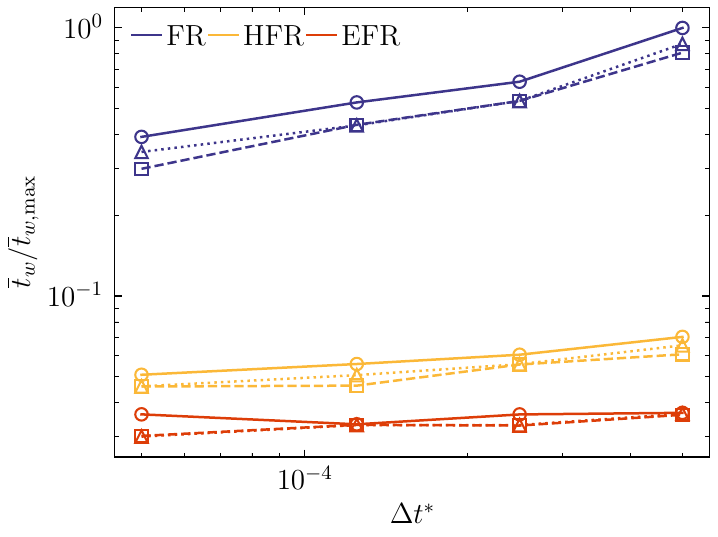}
  \end{subfigure}
  \caption{Normalized time to solve the linear system as a function of the time step size. The performance of FR methods has a larger dependence on the size of the time step size as opposed to hybridized methods. Markers for $c_{DG}$ (\usebox0), $c_{SD}$ (\usebox1), $c_{HU}$ (\usebox2) have been added at polynomial degrees $p=1$ to $p=4$}
  \label{fig:uadf-timestep}
\end{figure}

In terms of performance, Figure~\ref{fig:uadf-metrics} shows the error in the $L_2$ norm resulting from the $40\times40$ grids with all polynomial degrees against the wall-clock time. We observe the wall-clock time required to obtain a certain level of the $L_2$ error. Two values of the time step size were chosen to perform the comparison, namely \textcolor{\colornameedits}{ $\Delta t^* = 2.5\times 10^{-4}$ and $\Delta t^* = 5\times 10^{-5}$}, for which results appear in Figures~\ref{fig:uadf-metricsa} and \ref{fig:uadf-metricsb}, respectively. Clearly, the use of a large time step size is detrimental to the accuracy of the solution for the finest simulations in all runs. However, hybridized methods were able to achieve smaller levels of the $L_2$ norm in both cases at a fraction of the cost. Specifically, at $p=4$, simulations were between 14.06 and 17.52 times faster for HFR and between 22.5 to 30.1 times faster for EFR, compared to standard implicit FR schemes for the considered values of $c$. The speedups for all other runs on the finest grids are shown in Table~\ref{tab:speedupuadf}. Hence, reducing the time step size from $\Delta t^*=2.5\times10^{-4}$ to $\Delta t^* = 5\times 10^{-5}$, improved the performance of the standard FR formulation by $\sim1.6$ times per linear solve, but the improvements for the hybridized methods were not significant.
To further investigate the influence of the time step size on the performance results, we carry out an additional set of simulations on the $40\times40$ grid for a larger range of time step sizes at $p=3$. We show the results in Figure~\ref{fig:uadf-timestep}. Consistently, reducing the time step size improves the performance of the standard FR formulation per linear solve, where a reduction by a factor of ten showed an effect of 2.5 times better performance per step. However, the impact of the time step size was less significant for the hybridized methods, whose performance per linear solve remains almost constant with the considered $\Delta t$ increase, especially for the EFR schemes.
\begin{table}[htb]
\centering
\caption{Speedup factors for the unsteady linear advection-diffusion problem with 1600 elements. Values are shown for two time step sizes and compared against the corresponding FR formulation with the same $c$-parameter.}
\label{tab:speedupuadf}
\begin{tabular}{lllllllll}
\hline
                                   &     & \multicolumn{3}{l}{HFR}        &  & \multicolumn{3}{l}{EFR}        \\ \cline{3-5} \cline{7-9} 
$\Delta t^*$                         & $p$ & $c_{DG}$ & $c_{SD}$ & $c_{HU}$ &  & $c_{DG}$ & $c_{SD}$ & $c_{HU}$ \\ \hline
\multirow{4}{*}{$5\times 10^{-5}$} & 1   & 0.76     & 0.89     & 0.84     &  & 2.00     & 2.07     & 1.81     \\
                                   & 2   & 2.71     & 2.75     & 2.59     &  & 4.46     & 4.35     & 4.73     \\
                                   & 3   & 7.72     & 7.49     & 6.49     &  & 10.84    & 11.53    & 9.92     \\
                                   & 4   & 15.22    & 14.37    & 12.81    &  & 22.60    & 20.42    & 17.29    \\ \hline
\multirow{4}{*}{$2.5\times 10^{-4}$} & 1   & 1.08     & 1.08     & 1.03     &  & 2.95     & 2.58     & 2.24     \\
                                   & 2   & 4.41     & 4.25     & 3.52     &  & 7.78     & 6.78     & 6.00     \\
                                   & 3   & 10.43    & 9.59     & 9.62     &  & 17.42    & 16.18    & 16.26    \\
                                   & 4   & 17.52    & 16.51    & 14.06    &  & 30.12    & 26.35    & 22.56    \\ \hline
\end{tabular}
\end{table}

\subsection{Advection-Diffusion of a Gaussian Profile}
We now consider unsteady advection-diffusion of a Gaussian profile to analyze the stabilization mechanism of hybrid FR methods. To this end, we make use of a $[-5,5]^2$ domain with periodic boundary conditions and an initial condition
\begin{equation}
  u(\bm x, 0) = e^{-(x^2+y^2)}.
\end{equation}
We consider a pure advection problem with $\beta=0$, $\bm\alpha=[1,1]$ and an advection-diffusion problem with $\beta=0.01$ and the same advection velocity, \textcolor{\colornameedits}{following~\cite{sheshadriAnalysisStabilityFlux2018}}. We run a set of simulations using a third-order singly-diagonal Runge-Kutta (SDIRK) method with a small time step size $\Delta t^*=0.05$ on a $20\times20$ grid, representing a CFL number of 0.01. This small value helps mitigate the temporal errors. All runs are performed with $p=3$ spatial discretizations for five convective times $t^*=5$ for the advection problem and for one convective time $t^*=1$ for the advection-diffusion case.

\begin{figure}[htb]
  \centering
  \sbox0{\hwplotDG}\sbox1{\hwplotSD}\sbox2{\hwplotHU}
  \begin{subfigure}[b]{0.49\textwidth}
    \includegraphics[width=\textwidth]{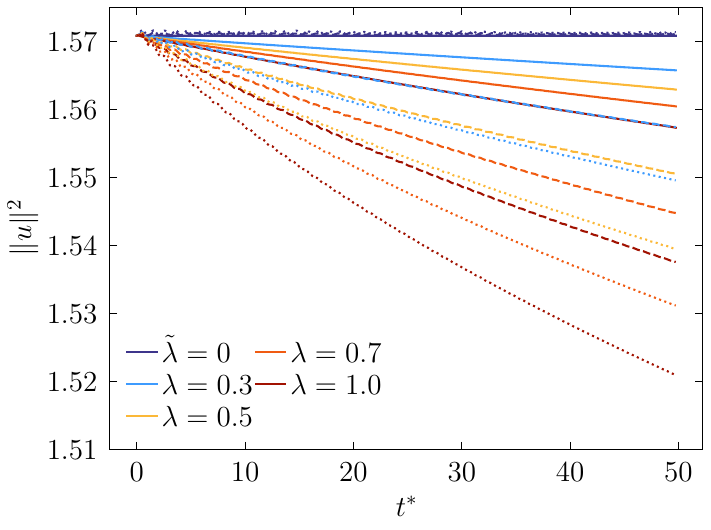}
    \caption{HFR}
  \end{subfigure}
    \begin{subfigure}[b]{0.49\textwidth}
    \includegraphics[width=\textwidth]{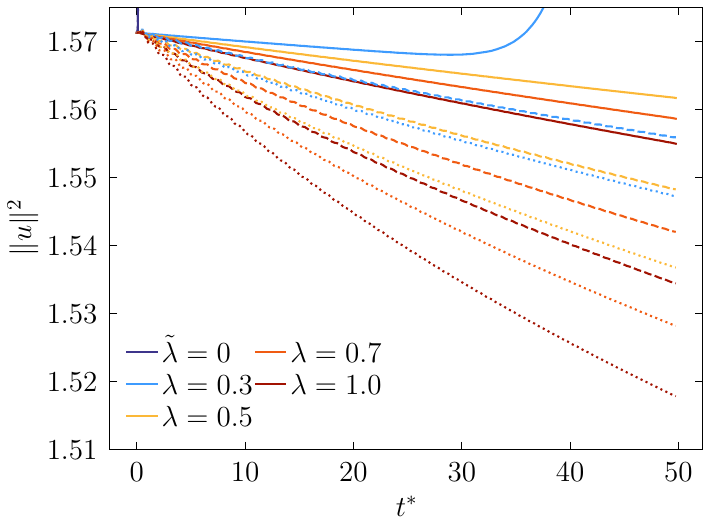}
    \caption{EFR}
  \end{subfigure}
  \caption{Time-evolution of the solution energy for hybridized FR methods for the Gaussian pulse case with $\beta=0$ for different values of $\lambda$. FR and HFR methods exhibit exactly the same curves so only one is shown on the left. The value $\tilde\lambda$ refers to the central FR method recovered via HFR with stabilization in Equations~\eqref{eq:hfrcentral1}-\eqref{eq:hfrcentral2} and $\lambda\rightarrow 0$ for EFR. Line strokes define methods with $c_{DG}$ (\usebox0), $c_{SD}$ (\usebox1), $c_{HU}$ (\usebox2)}
  \label{fig:stability}
\end{figure}

\begin{figure}[htb]
  \centering
  \sbox0{\hwplotDG}\sbox1{\hwplotSD}\sbox2{\hwplotHU}
  \begin{subfigure}[b]{0.49\textwidth}
    \includegraphics[width=\textwidth]{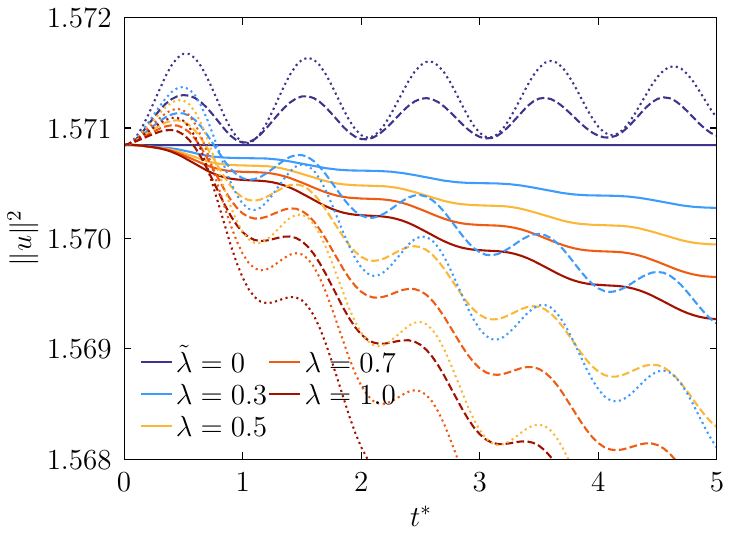}
    \caption{HFR}
  \end{subfigure}
    \begin{subfigure}[b]{0.49\textwidth}
    \includegraphics[width=\textwidth]{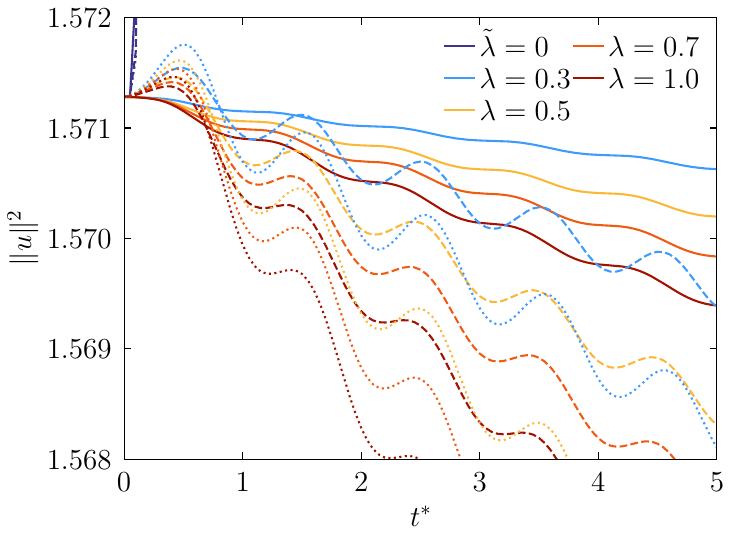}
    \caption{EFR}
  \end{subfigure}
  \caption{\textcolor{\colornameedits}{Zoom of the time-evolution of the solution energy for hybridized FR methods for the Gaussian pulse case with $\beta=0$ for different values of $\lambda$. The value $\tilde\lambda$ refers to the central FR method recovered via HFR with stabilization in Equations~\eqref{eq:hfrcentral1}-\eqref{eq:hfrcentral2} and $\lambda\rightarrow 0$ for EFR. Line strokes define methods with $c_{DG}$ (\usebox0), $c_{SD}$ (\usebox1), $c_{HU}$ (\usebox2)}}
  \label{fig:stabilityzoomed}
\end{figure}

First, we discuss the purely advective case. Recall that $\lambda=0$ results in undefined hybridized methods. Hence, for HFR methods, we have previously shown their equivalence in the stability analysis section when the conservation law includes only the advection operator. For the HFR method, we consider schemes ranging from the central approach of Equations~\eqref{eq:hfrcentral1}-\eqref{eq:hfrcentral2}, to increasing values of the upwinding parameter $\lambda\in\{0.3,0.5,0.7,1.0\}$, where $\lambda=1$ represents the upwind scheme. For the EFR method, we consider $\lambda\in\{0.001,0.3,0.5,0.7,1.0\}$. For all methods, we compute the evolution of the solution energy for correction parameters $c_{DG},c_{SD}$ and $c_{HU}$. Results are shown in Figure~\ref{fig:stability} with a zoomed-in version in Figure~\ref{fig:stabilityzoomed}. The behaviour in these figures is consistent with that previously observed in~\cite{vincentNewClassHighOrder2011,sheshadriStabilityFluxReconstruction2016}, where larger values of $c$ introduce additional dissipation and $\lambda\rightarrow0$ approaches a central FR method, for which the expected oscillatory behaviour of $c\neq0$ is observed. As opposed to HFR, EFR methods require a larger value of $\lambda$ for stability. As we previously discussed in the analysis of the stability section, we cannot easily guarantee stable EFR methods for $\lambda<\frac{1}{2}$ and, in fact, our experiments reveal that the minimum value is problem-dependent. Clearly, results for $\lambda=0.3$ are unstable for $c=0$, blowing up after about 35 convective times, but the added dissipation of $c_{SD}$ and $c_{HU}$ kept the simulations stable for this $\lambda$ parameter for the duration of these simulations. This is consistent with the analysis, where larger values of $c$ result \textcolor{\colornameedits}{in more negative $\frac{d}{dt}\Vert u\Vert_2$, i.e. more energy dissipation}. While EFR methods may appear to be stable for values $\lambda<\frac{1}{2}$ for some problems, they may be only mildly stable and may blow up later in long-time integration simulations. 
\begin{figure}[htb]
  \centering
  \sbox0{\hwplotDG}\sbox1{\hwplotSD}\sbox2{\hwplotHU}
  \begin{subfigure}[b]{0.49\textwidth}
    \includegraphics[width=\textwidth]{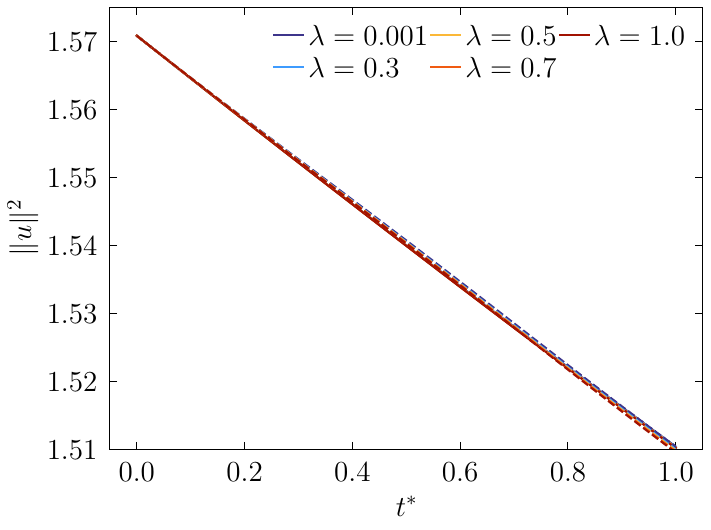}
    \caption{FR}
  \end{subfigure}
    \begin{subfigure}[b]{0.49\textwidth}
    \includegraphics[width=\textwidth]{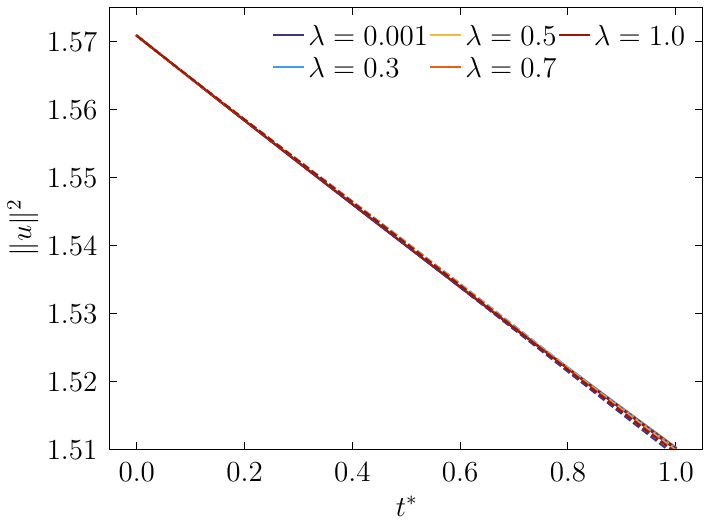}
    \caption{HFR}
  \end{subfigure}
      \begin{subfigure}[b]{0.49\textwidth}
    \includegraphics[width=\textwidth]{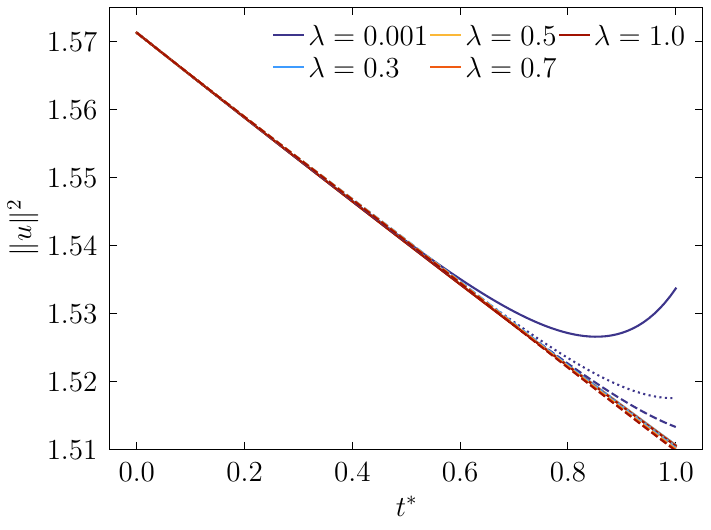}
    \caption{EFR}
  \end{subfigure}
  \caption{Time-evolution of the solution energy for hybridized FR methods for the Gaussian pulse case with $\beta=0.01$ for different values of $\lambda$. Line strokes define methods with $c_{DG}$ (\usebox0), $c_{SD}$ (\usebox1), $c_{HU}$ (\usebox2)}
  \label{fig:advdiffstability}
\end{figure}
\begin{figure}[htb]
  \centering
  \sbox0{\hwplotDG}\sbox1{\hwplotSD}\sbox2{\hwplotHU}
      \begin{subfigure}[b]{0.49\textwidth}
    \includegraphics[width=\textwidth]{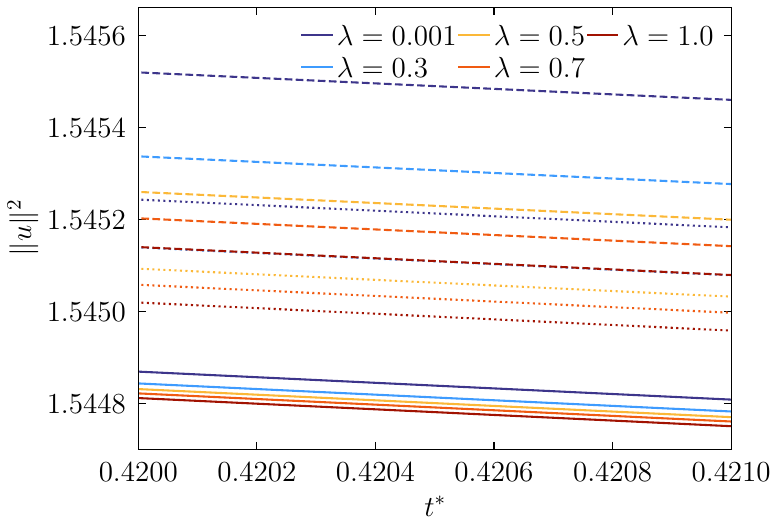}
    \caption{FR}
  \end{subfigure}
  \begin{subfigure}[b]{0.49\textwidth}
    \includegraphics[width=\textwidth]{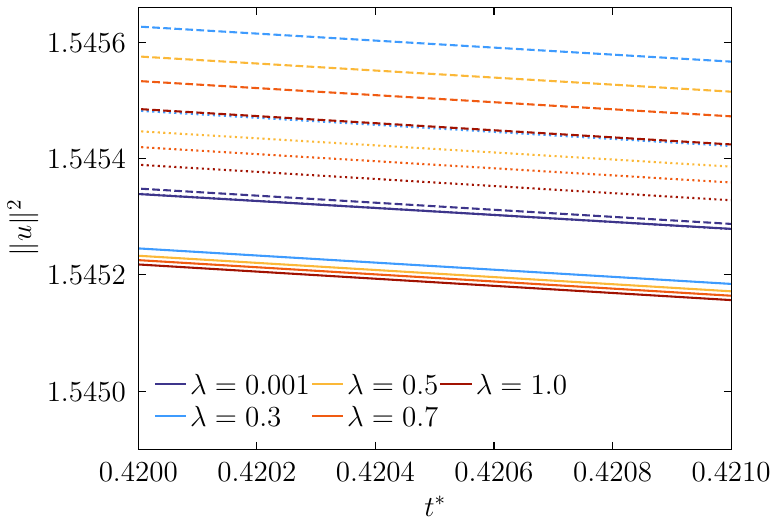}
    \caption{HFR}
  \end{subfigure}
    \begin{subfigure}[b]{0.49\textwidth}
    \includegraphics[width=\textwidth]{stability_advdiff_zoom_efr}
    \caption{EFR}{}
  \end{subfigure}
  \caption{Zoom of the time-evolution of the solution energy for hybridized FR methods for the Gaussian pulse case with $\beta=0.01$ for different values of $\lambda$. Line strokes define methods with $c_{DG}$ (\usebox0), $c_{SD}$ (\usebox1), $c_{HU}$ (\usebox2)}
  \label{fig:advdiffstabilityzoomed}
\end{figure}

Next, we consider the advection-diffusion case on the same computational grid. Similarly, we plot the energy of the solution against the convective time in Figure~\ref{fig:advdiffstability}. In this case, we show results for FR, HFR and EFR, none of which are equivalent. In the case of FR, larger $\lambda$ means larger numerical dissipation, consistent with the purely advective case. Due to the physical diffusion, the energy of the solution is seen to decay for the linear advection-diffusion problem with $\beta=0.01$. A larger range of $\lambda$ parameters is stable for the EFR method compared to the purely advective case due to physical dissipation. However, it can be seen that for the smallest considered value, an increase in energy is seen close to the end of the simulation for all considered values of $c$, with larger $c$ introducing additional dissipation. Results look very similar among FR, HFR and EFR from this view, and hence we show a zoomed-in version for a time range $0.4 < t^* < 0.45$ in Figure~\ref{fig:advdiffstabilityzoomed}, moments before the peak of the Gaussian pulse reaches the periodic boundary. Smaller values of $\lambda\rightarrow0$ for HFR and EFR  make the gradient jump term grow in the trace equation. See Equations~\eqref{eq:explicittracehfr}, and \eqref{eq:explicittraceefr}. Then, setting $\lambda\rightarrow0$ does not recover the behavior of a central scheme when diffusion operators appear. Hence, for advection-diffusion problems, the HFR central scheme cannot be recovered with the formulation in~\eqref{eq:hfrcentral1}, \eqref{eq:hfrcentral2} for the advection-diffusion case. With this example, we observe that results are consistent with the analytical findings of the stability section and demonstrate a range of stable HFR and EFR methods for advection-diffusion. 

\subsection{Planar Couette flow}
\begin{figure}[htb]
\sbox0{\hwplotDGmarkerfilled}\sbox1{\hwplotSDmarkerfilled}\sbox2{\hwplotHUmarkerfilled}%
\centering
    \begin{subfigure}[b]{0.49\textwidth}
    \includegraphics[width=\textwidth]{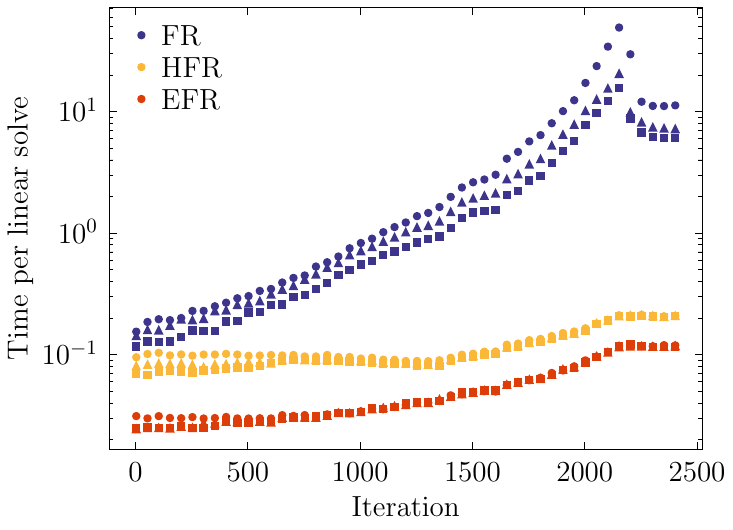}
      \label{fig:couette-tksp}
    \end{subfigure}
  \caption{Average time to solve the linear system per iteration ($p=4$ finest grids). Markers represent methods with $c_{DG}$ (\usebox0), $c_{SD}$ (\usebox1), $c_{HU}$ (\usebox2)}
  \label{fig:couette-cumulative-p4}
\end{figure}

Finally, we present a problem involving the compressible Navier-Stokes equations. Planar Couette flow is a well-known two-dimensional case to perform verification of the viscous fluxes given its simplification of the Navier-Stokes equations. This problem consists of viscous flow between two plates separated by a distance of $\delta$. A moving wall is located at $y=\delta$ with temperature $T_e$ and constant velocity $v_e$, which drives the flow in the positive $x$-direction. At $y = 0$, a fixed wall $(v_w = 0)$ with temperature $T_w$ is placed. Due to the no-slip condition, the flow variables are equal to those of the walls at $y = 0$ and $y = \delta$, respectively. Hence the flow experiences a temperature gradient due to viscous dissipation. The exact temperature profile can be computed from\textcolor{\colornameedits}{~\cite{andersonFundamentalsAerodynamics2010,castonguayEnergyStableFlux2013}}
\begin{equation}
  T = T_w + \left[T_e - T_w + \frac{\operatorname{Pr}}{2 c_p} v_e^2\right] \frac{y}{\delta} - \frac{\operatorname{Pr}}{2 c_p} v_e^2 \left(\frac{y}{\delta}\right)^2.
\end{equation}
where $\operatorname{Pr}=0.71$ \textcolor{\colornameedits}{and $c_p$ is the specific heat at constant pressure corresponding to specific heat ratio $\gamma=1.4$}. The Mach number is set to 0.1 and $\operatorname{Re}=5$.
\begin{figure}[htb]
  \centering
  \begin{subfigure}[b]{0.49\textwidth}
    \includegraphics[width=\textwidth]{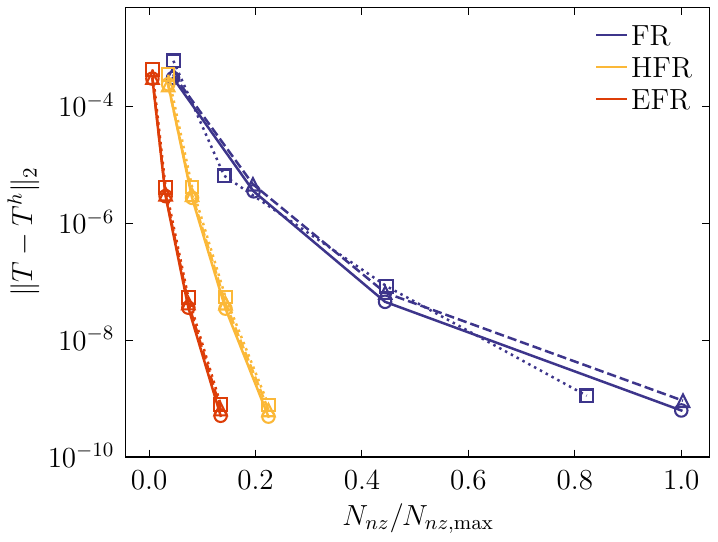}
    \caption{Relative number of nonzeros}
    \label{fig:couette-metricsa}
  \end{subfigure}
    \begin{subfigure}[b]{0.49\textwidth}
    \includegraphics[width=\textwidth]{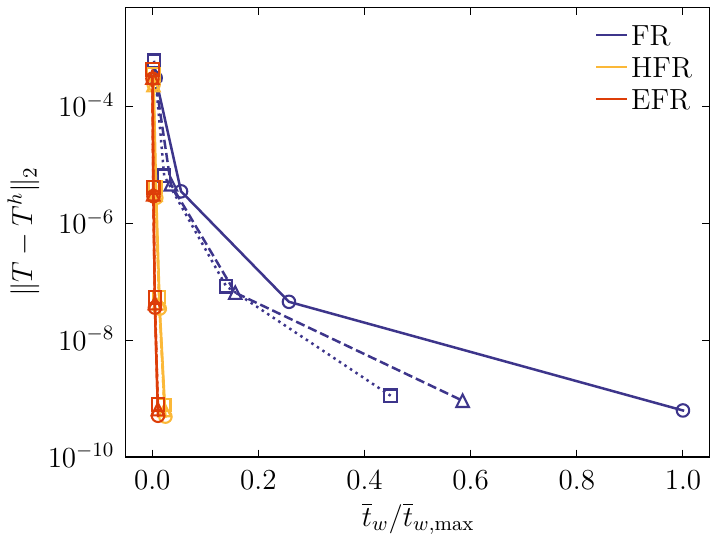}
    \caption{Relative solver time}
    \label{fig:couette-metricsb}
  \end{subfigure}
  \caption{Performance metrics for the planar Couette problem}
  \label{fig:couette-metrics}
\end{figure}
We consider four levels of refinement using grids composed of $4\times2$, $8\times 4$, $16\times 8$ and $32\times 16$ quadrilateral elements with $p=1-4$ schemes and the same $c$-parameters considered in the previous section. For stabilization, we consider a Lax-Friedrichs type matrix for the inviscid fluxes such as that in~\cite{pereiraPerformanceAccuracyHybridized2022} for the Euler equations, and constant viscous stabilization of the form
\begin{equation}
  s_{f,k}^{(v)} = \frac{1}{\operatorname{Re}},
\end{equation}
for the hybridized formulations. For standard FR, we also make use of a Lax-Friedrichs Riemann solver for the inviscid fluxes and the LDG method with an arbitrary directional switch for the viscous component. The simulations were initialized with a stationary problem $(\bm v=\bm 0)$ and allowed to run until the density residual converged to a tolerance of $10^{-10}$ in the $L_\infty$ norm, which proved to be sufficient for the temperature error in the $L_2$ norm to converge. For convergence acceleration, a relaxation factor was used with an implicit Euler scheme to converge the nonlinear residuals via a time step ramp function. Specifically, the following function at the $i$-th iteration was used to update the time step size
\begin{equation}
\Delta t_i = \begin{cases}
\Delta t_0  & \quad i<20, \\
2^{1/16}\Delta t_{i-1} & \mod(i,5)=0~\text{and}~ \Delta t_{i-1} < 10^4 \Delta t_0, \\
\Delta t_{i-1} & \quad \text{otherwise},
\end{cases}
\end{equation}
\begin{table}[ht]
\centering
\caption{Speedup factors for the Couette problem with $32\times16$ elements for the HFR and EFR methods. These factors take into account the time to solve the global and local problems, but exclude the assembly time of the Jacobian matrices}
\label{tab:speedupcouette}
\begin{tabular}{llllllll}
\hline
    & \multicolumn{3}{l}{HFR}                                                                    &  & \multicolumn{3}{l}{EFR}                                                                    \\ \cline{2-4} \cline{6-8} 
$p$ & \multicolumn{1}{l}{$c_{DG}$} & \multicolumn{1}{l}{$c_{SD}$} & \multicolumn{1}{l}{$c_{HU}$} &  & \multicolumn{1}{l}{$c_{DG}$} & \multicolumn{1}{l}{$c_{SD}$} & \multicolumn{1}{l}{$c_{HU}$} \\ \hline
1   & 2.31                         & 1.72                         & 1.53                         &  & 14.84                        & 11.11                         & 9.33                        \\
2   & 6.91                         & 5.46                         & 3.57                         &  & 21.95                       & 16.66                         & 10.74                         \\
3   & 18.10                         & 11.73                         & 10.77                         &  & 49.45                       & 30.31                       & 27.26                        \\
4   & 40.38                        & 25.12                        & 19.86                        &  & 90.7                        & 54.24                        & 41.73                       \\ \hline
\end{tabular}
\end{table}which was chosen empirically to accelerate convergence of the FR simulation and was used for all grids and polynomial degrees. The value of the base time step size was set to \textcolor{\colornameedits}{$\Delta t_0^*=\Delta t/t_c=1.1832\times10^{-6}$, with characteristic time $t_c=\delta/U_\infty$, where $U_\infty$ is the freestream velocity based on Mach}. The Jacobian matrix was computed exactly and updated every five time steps to reduce the computational cost associated with its assembly. Verification is presented in Table~\ref{tab:ordercouette}, where the $L_2$ norm of the error is computed for the aforementioned levels of refinement. It is interesting to see that for this nonlinear problem, HFR displayed the lowest $L_2$-norm, followed by EFR and then FR. Specifically, $L_2$ errors of $6.17\times10^{-10},~4.75\times10^{-10}$ and $4.86\times10^{-10}$ were obtained for FR, HFR and EFR, respectively. The expected orders of accuracy were obtained in all cases and considered values of $c$. Consistent with our previous experiments, we analyze the performance based on the time spent on the solution of the linear system, and for the hybridized method, this accounts for the solution of the local problems. The number of nonzeros in the implicit system reduces between 4 to 7 times for the HFR and EFR methods compared to FR for the finest problems, as shown in Figure~\ref{fig:couette-metricsa}. Interestingly, the effect of the $c$-parameter on the FR simulations can be significant. For $c=c_{DG}$, simulations are twice as expensive as those using $c_{HU}$. This can be attributed to fewer nonzero entries for this value of $c$ as well as a stiffer problem resulting from setting $c=0$. This is consistent with our observations in~\cite{pereiraPerformanceAccuracyHybridized2022} for in the advection regime. We observed significant speedup values, which are shown in Table~\ref{tab:speedupcouette}. For example, between 19.86 and 40.38 times faster simulations were observed for the HFR problems and between 41.73 to 90.7 times faster solutions for the EFR problem. 
\begin{figure}[htb]
  \sbox0{\hwplotDG}\sbox1{\hwplotSD}\sbox2{\hwplotHU}
  \centering
  \begin{subfigure}[t]{0.48\textwidth}
    \includegraphics[width=\textwidth]{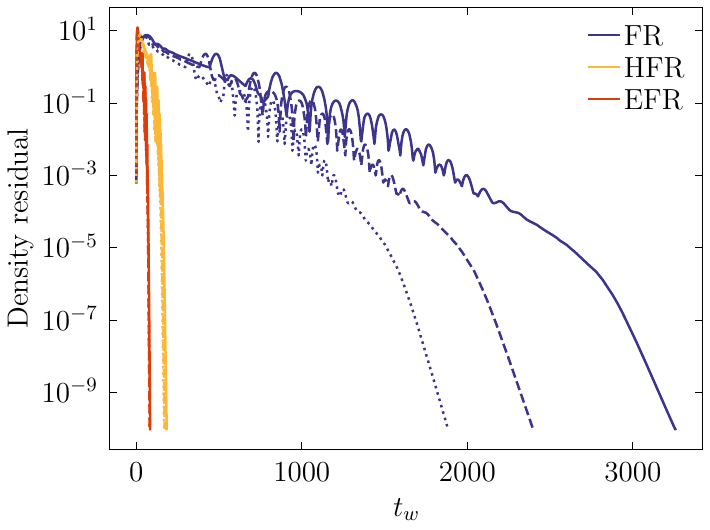}
    \caption{Density residual}
      \label{fig:couette-residual}
  \end{subfigure}
    \begin{subfigure}[t]{0.49\textwidth}
    \includegraphics[width=\textwidth]{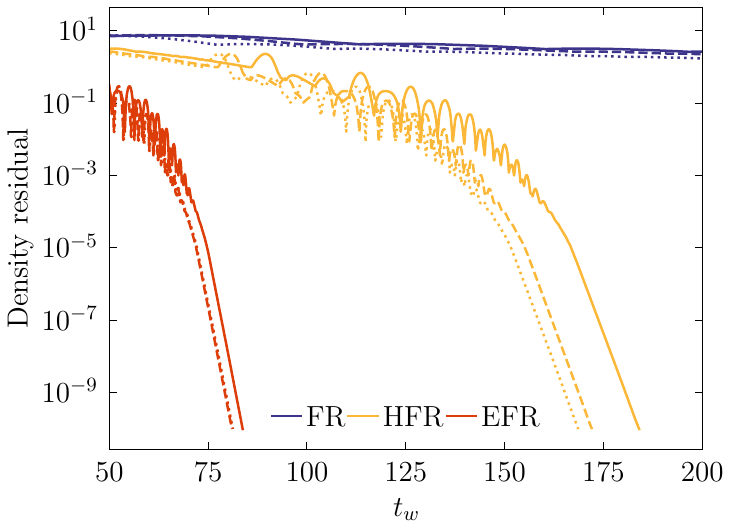}
    \caption{\textcolor{\colornameedits}{Density residual zoom}}
      \label{fig:couette-residual-zoom}
  \end{subfigure}

    \begin{subfigure}[t]{0.48\textwidth}
    \includegraphics[width=\textwidth]{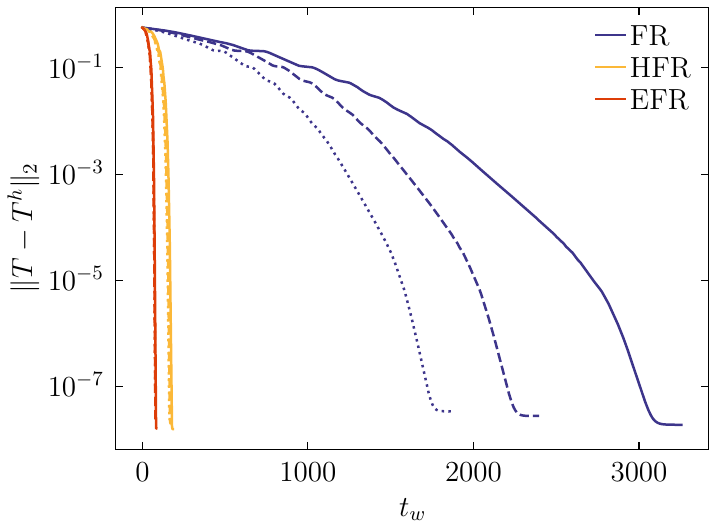}
    \caption{Temperature $L_2$ error}
      \label{fig:couette-error}
    \end{subfigure}
        \begin{subfigure}[t]{0.49\textwidth}
    \includegraphics[width=\textwidth]{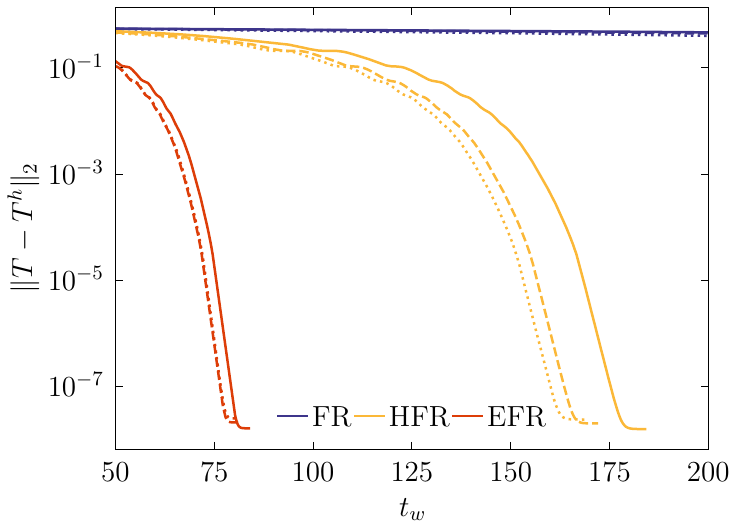}
    \caption{\textcolor{\colornameedits}{Temperature $L_2$ error zoom}}
      \label{fig:couette-error-zoom}
    \end{subfigure}
  \caption{Temperature and density residual against wall-clock time in seconds for a $p=4$ scheme with different correction functions on the $16\times8$ grid. Line strokes represent methods with $c_{DG}$ (\usebox0), $c_{SD}$ (\usebox1), $c_{HU}$ (\usebox2)}
  \label{fig:couette-p4}
\end{figure}
Plots of the residual and temperature error against the wall-clock time are shown in Figure~\ref{fig:couette-p4} for the finest grids and $p=4$ runs as an example of the evolution of the convergence. The speedups obtained with this problem are significantly larger than the ratios of nonzeros. In Figure~\ref{fig:couette-cumulative-p4}, we show the time spent to solve the linear system at every iteration. Here, we observe that the increasing time step size was detrimental to the FR time spent on solving these systems. Much more so than for the hybridized methods. Hence, we can significantly reduce the cost of implicit FR simulations for nonlinear viscous problems via hybridization.
\section{Conclusions}\label{sec:conclusion}
In this work, we have studied the performance, accuracy and stability of hybridized flux reconstruction methods for advection-diffusion problems. First, we discussed implementation details for hybridization of these methods. We then demonstrated linear stability constraints for advection-diffusion problems. We found that for pure advection, HFR methods can fully recover standard FR formulations even for the central case, where a special form of stabilization is required. For advection-diffusion problems, HFR methods do not recover FR-LDG approaches for finite stabilization since the common value of the solution contains the jump of the auxiliary variable. It was also observed that the minimum stabilization for EFR methods depends on the problem and the value of $c$, since the latter acts as an added dissipation mechanism. For a series of numerical experiments, hybridized methods displayed superior accuracy and computational performance. The considered schemes achieved up to 30 times faster solutions in unsteady problems and up to 90 times faster for the steady nonlinear planar Couette case. Since hybridized FR methods for diffusion present a superconvergence property, we redefined the post-processing schemes for consistency with FR formulations by introducing correction functions. It was found that HFR methods  recover the super accuracy property of HDG methods, but they do not possess superconvergence if $c\neq c_{DG}$ and $p=1$. \textcolor{\colornameedits}{Further analysis on the superconvergence properties under more general correction fields and the extension of the stability analysis to curved discretizations on more general element types are subjects of ongoing work.}

\section*{Acknowledgements}
\noindent We acknowledge the support of the Natural Sciences and Engineering Research Council of Canada (NSERC), [RGPAS-2017-507988, RGPIN-2017-06773], Concordia University via the Team Seed pro- gram and the Fonds de Recherche du Quebec - Nature et Technologie (FRQNT) via a B2X scholarship. This research was enabled in part by support provided the Digital Research Alliance of Canada (www.alliancecan.ca) via a Resources for Research Groups allocation.

\section*{Data Statement}
\noindent Data relating to the results in this manuscript can be downloaded from the publication's website under a CC BY-NC-ND 4.0 license.

\pagebreak

\begin{sidewaystable}[h]
\appendix
\setcounter{table}{0}
\section{Grid Convergence Tables}
\tiny
\centering
\setlength\tabcolsep{1pt}
\renewcommand{\arraystretch}{1.2}
\caption{$L_2$ norm of the error for the steady-state linear advection-diffusion problem}
\label{tab:ordersadvdiff}
\begin{tabular}{llllllllllllllllllllll}
\hline
                   &              & \multicolumn{6}{l}{FR}                                 &  & \multicolumn{6}{l}{HFR}                                &  & \multicolumn{6}{l}{EFR}                                \\ \cline{3-8} \cline{10-15} \cline{17-22} 
$p$                & Mesh size    & $c_{DG}$ & Order & $c_{SD}$ & Order & $c_{HU}$ & Order &  & $c_{DG}$ & Order & $c_{SD}$ & Order & $c_{HU}$ & Order &  & $c_{DG}$ & Order & $c_{SD}$ & Order & $c_{HU}$ & Order \\ \hline
\multirow{4}{*}{1} & 5$\times$5   & 5.79E-03 & -     & 6.39E-03 & -     & 9.82E-03 & -     &  & 4.02E-03 & -     & 4.94E-03 & -     & 7.92E-03 & -     &  & 6.49E-03 & -     & 7.15E-03 & -     & 7.87E-03 & -     \\
                   & 10$\times$10 & 1.31E-03 & 2.14  & 1.57E-03 & 2.03  & 2.62E-03 & 1.91  &  & 1.09E-03 & 1.88  & 1.23E-03 & 2.01  & 2.07E-03 & 1.93  &  & 1.45E-03 & 2.16  & 1.56E-03 & 2.20  & 1.64E-03 & 2.26  \\
                   & 20$\times$20 & 3.01E-04 & 2.12  & 3.83E-04 & 2.03  & 6.98E-04 & 1.91  &  & 2.25E-04 & 2.27  & 3.11E-04 & 1.98  & 6.58E-04 & 1.65  &  & 4.94E-04 & 1.55  & 4.99E-04 & 1.65  & 4.74E-04 & 1.79  \\
                   & 40$\times$40 & 7.78E-05 & 1.95  & 9.58E-05 & 2.00  & 1.75E-04 & 2.00  &  & 6.80E-05 & 1.73  & 9.48E-05 & 1.71  & 2.03E-04 & 1.69  &  & 1.67E-04 & 1.57  & 1.59E-04 & 1.65  & 1.34E-04 & 1.82  \\ \hline
\multirow{4}{*}{2} & 5$\times$5   & 1.37E-03 & -     & 1.29E-03 & -     & 1.33E-03 & -     &  & 6.94E-04 & -     & 5.95E-04 & -     & 8.84E-04 & -     &  & 7.78E-04 & -     & 5.54E-04 & -     & 5.73E-04 & -     \\
                   & 10$\times$10 & 1.54E-04 & 3.15  & 1.71E-04 & 2.92  & 2.06E-04 & 2.69  &  & 1.10E-04 & 2.66  & 1.29E-04 & 2.20  & 1.80E-04 & 2.30  &  & 1.09E-04 & 2.84  & 1.04E-04 & 2.41  & 1.20E-04 & 2.25  \\
                   & 20$\times$20 & 1.98E-05 & 2.96  & 2.32E-05 & 2.88  & 2.98E-05 & 2.79  &  & 1.66E-05 & 2.73  & 2.30E-05 & 2.49  & 3.44E-05 & 2.39  &  & 1.55E-05 & 2.81  & 1.78E-05 & 2.55  & 2.36E-05 & 2.35  \\
                   & 40$\times$40 & 2.42E-06 & 3.03  & 2.96E-06 & 2.97  & 4.03E-06 & 2.89  &  & 2.30E-06 & 2.85  & 3.42E-06 & 2.75  & 5.27E-06 & 2.70  &  & 2.09E-06 & 2.89  & 2.65E-06 & 2.75  & 3.70E-06 & 2.67  \\ \hline
\multirow{4}{*}{3} & 5$\times$5   & 2.43E-04 & -     & 2.48E-04 & -     & 2.62E-04 & -     &  & 1.69E-04 & -     & 1.63E-04 & -     & 1.90E-04 & -     &  & 1.72E-04 & -     & 1.50E-04 & -     & 1.59E-04 & -     \\
                   & 10$\times$10 & 1.60E-05 & 3.92  & 1.93E-05 & 3.68  & 2.27E-05 & 3.53  &  & 1.29E-05 & 3.71  & 1.73E-05 & 3.24  & 2.20E-05 & 3.11  &  & 1.24E-05 & 3.79  & 1.47E-05 & 3.35  & 1.78E-05 & 3.16  \\
                   & 20$\times$20 & 1.01E-06 & 3.99  & 1.28E-06 & 3.91  & 1.56E-06 & 3.86  &  & 9.09E-07 & 3.83  & 1.34E-06 & 3.69  & 1.77E-06 & 3.64  &  & 8.21E-07 & 3.92  & 1.06E-06 & 3.80  & 1.35E-06 & 3.72  \\
                   & 40$\times$40 & 6.16E-08 & 4.03  & 8.10E-08 & 3.98  & 1.02E-07 & 3.94  &  & 5.99E-08 & 3.92  & 9.37E-08 & 3.84  & 1.26E-07 & 3.81  &  & 5.01E-08 & 4.03  & 6.77E-08 & 3.97  & 8.95E-08 & 3.91  \\ \hline
\multirow{4}{*}{4} & 5$\times$5   & 3.69E-05 & -     & 4.11E-05 & -     & 4.46E-05 & -     &  & 2.62E-05 & -     & 3.26E-05 & -     & 3.87E-05 & -     &  & 2.59E-05 & -     & 2.89E-05 & -     & 3.30E-05 & -     \\
                   & 10$\times$10 & 1.37E-06 & 4.75  & 1.69E-06 & 4.61  & 1.92E-06 & 4.54  &  & 1.18E-06 & 4.48  & 1.68E-06 & 4.28  & 2.03E-06 & 4.25  &  & 1.14E-06 & 4.51  & 1.48E-06 & 4.28  & 1.73E-06 & 4.25  \\
                   & 20$\times$20 & 4.33E-08 & 4.98  & 5.73E-08 & 4.88  & 6.74E-08 & 4.83  &  & 4.10E-08 & 4.85  & 6.29E-08 & 4.74  & 7.81E-08 & 4.70  &  & 3.90E-08 & 4.87  & 5.48E-08 & 4.76  & 6.60E-08 & 4.72  \\
                   & 40$\times$40 & 1.34E-09 & 5.01  & 1.79E-09 & 5.00  & 2.14E-09 & 4.98  &  & 1.35E-09 & 4.93  & 2.12E-09 & 4.89  & 2.66E-09 & 4.88  &  & 1.27E-09 & 4.95  & 1.83E-09 & 4.90  & 2.24E-09 & 4.88  \\ \hline
\end{tabular}

\caption{$L_2$-norm of the post-processed solution error for the steady-state linear diffusion problem}
\label{tab:orderzsadf}
\begin{tabular}{lllllllllllllll}
\hline
                   &              & \multicolumn{6}{l}{FR}                                 &  & \multicolumn{6}{l}{HFR}                                \\ \cline{3-8} \cline{10-15} 
$p$                & Mesh size    & $c_{DG}$ & Order & $c_{SD}$ & Order & $c_{HU}$ & Order &  & $c_{DG}$ & Order & $c_{SD}$ & Order & $c_{HU}$ & Order \\ \hline
\multirow{4}{*}{1} & 10$\times$10 & 5.39E-03 & -     & 7.46E-03 & -     & 1.70E-02 & -     &  & 2.22E-03 & -     & 2.99E-03 & -     & 7.37E-03 & -     \\
                   & 20$\times$20 & 9.84E-04 & 2.14  & 1.34E-03 & 2.03  & 3.06E-03 & 1.91  &  & 2.39E-04 & 3.22  & 5.50E-04 & 2.44  & 1.86E-03 & 1.98  \\
                   & 40$\times$40 & 2.04E-04 & 2.12  & 3.08E-04 & 2.03  & 6.96E-04 & 1.91  &  & 3.23E-05 & 2.89  & 1.22E-04 & 2.18  & 4.71E-04 & 1.98  \\
                   & 80$\times$80 & 5.22E-05 & 1.95  & 7.23E-05 & 2.00  & 1.60E-04 & 2.00  &  & 4.33E-06 & 2.90  & 3.47E-05 & 1.81  & 1.38E-04 & 1.77  \\ \hline
\multirow{4}{*}{2} & 10$\times$10 & 1.84E-03 & -     & 1.90E-03 & -     & 1.99E-03 & -     &  & 3.57E-04 & -     & 3.73E-04 & -     & 4.98E-04 & -     \\
                   & 20$\times$20 & 1.21E-04 & 3.15  & 1.37E-04 & 2.92  & 1.56E-04 & 2.69  &  & 2.46E-05 & 3.86  & 3.13E-05 & 3.58  & 4.35E-05 & 3.52  \\
                   & 40$\times$40 & 1.19E-05 & 2.96  & 1.36E-05 & 2.88  & 1.57E-05 & 2.79  &  & 1.81E-06 & 3.77  & 2.61E-06 & 3.58  & 3.81E-06 & 3.52  \\
                   & 80$\times$80 & 1.36E-06 & 3.03  & 1.59E-06 & 2.97  & 1.87E-06 & 2.89  &  & 1.25E-07 & 3.86  & 1.90E-07 & 3.78  & 2.84E-07 & 3.75  \\ \hline
\multirow{4}{*}{3} & 10$\times$10 & 2.72E-04 & -     & 3.52E-04 & -     & 3.99E-04 & -     &  & 4.80E-05 & -     & 5.02E-05 & -     & 6.40E-05 & -     \\
                   & 20$\times$20 & 1.04E-05 & 3.92  & 1.39E-05 & 3.68  & 1.60E-05 & 3.53  &  & 1.77E-06 & 4.76  & 1.86E-06 & 4.75  & 2.52E-06 & 4.67  \\
                   & 40$\times$40 & 5.02E-07 & 3.99  & 6.75E-07 & 3.91  & 7.83E-07 & 3.86  &  & 5.20E-08 & 5.09  & 6.08E-08 & 4.94  & 9.66E-08 & 4.70  \\
                   & 80$\times$80 & 2.83E-08 & 4.03  & 3.78E-08 & 3.98  & 4.41E-08 & 3.94  &  & 1.56E-09 & 5.06  & 2.07E-09 & 4.87  & 3.55E-09 & 4.77  \\ \hline
\multirow{4}{*}{4} & 10$\times$10 & 3.53E-05 & -     & 4.28E-05 & -     & 4.61E-05 & -     &  & 7.73E-06 & -     & 9.81E-06 & -     & 1.13E-05 & -     \\
                   & 20$\times$20 & 8.25E-07 & 4.75  & 1.03E-06 & 4.61  & 1.13E-06 & 4.54  &  & 1.49E-07 & 5.70  & 2.13E-07 & 5.53  & 2.53E-07 & 5.48  \\
                   & 40$\times$40 & 2.05E-08 & 4.98  & 2.52E-08 & 4.88  & 2.76E-08 & 4.83  &  & 2.63E-09 & 5.82  & 3.99E-09 & 5.74  & 4.79E-09 & 5.72  \\
                   & 80$\times$80 & 6.27E-10 & 5.01  & 7.70E-10 & 5.00  & 8.44E-10 & 4.98  &  & 5.00E-11 & 5.72  & 7.01E-11 & 5.83  & 8.34E-11 & 5.85  \\ \hline
\end{tabular}
\end{sidewaystable}

\begin{sidewaystable}[htb]
\tiny
\centering
\setlength\tabcolsep{4pt}
\renewcommand{\arraystretch}{1.2}
\caption{$L_2$-norm of the solution error for the steady-state linear diffusion problem with consistent LDG switch}
\label{tab:orderzsadvdiffconsistent}
\begin{tabular}{lllllllllllllll}
\hline
                   &              & \multicolumn{6}{l}{FR (no post-processing)}             &  & \multicolumn{6}{l}{FR (post-processed)}                 \\ \cline{3-15} 
$p$                & Mesh size    & $c_{DG}$ & Order & $c_{SD}$ & Order & $c_{HU}$ & Order &  & $c_{DG}$ & Order & $c_{SD}$ & Order & $c_{HU}$ & Order \\ \cline{1-8} \cline{10-15} 
\multirow{4}{*}{1} & 5$\times$5   & 8.56E-03 & -     & 9.86E-03 & -     & 1.33E-02 & -     &  & 1.03E-01 & -     & 1.24E-01 & -     & 1.98E-01 & -     \\
                   & 10$\times$10 & 1.24E-03 & 2.79  & 1.75E-03 & 2.49  & 3.46E-03 & 1.95  &  & 2.86E-02 & 3.26  & 3.28E-02 & 2.22  & 5.47E-02 & 2.11  \\
                   & 20$\times$20 & 3.13E-04 & 1.98  & 4.30E-04 & 2.03  & 8.45E-04 & 2.03  &  & 5.93E-03 & 2.91  & 8.07E-03 & 2.44  & 1.66E-02 & 2.29  \\
                   & 40$\times$40 & 8.03E-05 & 1.96  & 1.12E-04 & 1.94  & 2.31E-04 & 1.87  &  & 1.78E-03 & 2.90  & 2.45E-03 & 1.98  & 5.09E-03 & 1.92  \\ \cline{1-8} \cline{10-15} 
\multirow{4}{*}{2} & 5$\times$5   & 1.89E-03 & -     & 2.14E-03 & -     & 2.44E-03 & -     &  & 1.74E-02 & -     & 1.69E-02 & -     & 2.64E-02 & -     \\
                   & 10$\times$10 & 1.64E-04 & 3.53  & 2.00E-04 & 3.42  & 2.55E-04 & 3.26  &  & 2.87E-03 & 4.14  & 3.47E-03 & 3.97  & 4.93E-03 & 3.83  \\
                   & 20$\times$20 & 2.06E-05 & 2.99  & 2.82E-05 & 2.83  & 4.02E-05 & 2.67  &  & 4.35E-04 & 4.10  & 6.10E-04 & 4.02  & 9.20E-04 & 3.90  \\
                   & 40$\times$40 & 2.57E-06 & 3.01  & 3.76E-06 & 2.91  & 5.65E-06 & 2.83  &  & 6.02E-05 & 4.07  & 9.00E-05 & 4.02  & 1.39E-04 & 3.94  \\ \cline{1-8} \cline{10-15} 
\multirow{4}{*}{3} & 5$\times$5   & 3.04E-04 & -     & 3.81E-04 & -     & 4.41E-04 & -     &  & 4.36E-03 & -     & 4.38E-03 & -     & 5.27E-03 & -     \\
                   & 10$\times$10 & 1.70E-05 & 4.16  & 2.25E-05 & 4.08  & 2.74E-05 & 4.00  &  & 3.37E-04 & 4.96  & 4.61E-04 & 4.92  & 5.93E-04 & 4.96  \\
                   & 20$\times$20 & 1.04E-06 & 4.03  & 1.51E-06 & 3.89  & 1.96E-06 & 3.81  &  & 2.38E-05 & 5.11  & 3.54E-05 & 5.21  & 4.70E-05 & 5.30  \\
                   & 40$\times$40 & 6.47E-08 & 4.01  & 9.93E-08 & 3.93  & 1.32E-07 & 3.89  &  & 1.57E-06 & 5.05  & 2.46E-06 & 5.11  & 3.33E-06 & 5.17  \\ \cline{1-8} \cline{10-15} 
\multirow{4}{*}{4} & 5$\times$5   & 4.82E-05 & -     & 6.22E-05 & -     & 7.02E-05 & -     &  & 6.87E-04 & -     & 8.92E-04 & -     & 1.07E-03 & -     \\
                   & 10$\times$10 & 1.43E-06 & 5.08  & 1.98E-06 & 4.98  & 2.33E-06 & 4.91  &  & 3.08E-05 & 5.97  & 4.47E-05 & 6.04  & 5.44E-05 & 6.09  \\
                   & 20$\times$20 & 4.49E-08 & 4.99  & 6.72E-08 & 4.88  & 8.23E-08 & 4.82  &  & 1.07E-06 & 6.04  & 1.66E-06 & 6.12  & 2.06E-06 & 6.18  \\
                   & 40$\times$40 & 1.45E-09 & 4.95  & 2.21E-09 & 4.93  & 2.74E-09 & 4.91  &  & 3.44E-08 & 5.35  & 5.53E-08 & 5.56  & 6.97E-08 & 5.79  \\ \cline{1-8} \cline{10-15} 
\end{tabular}
\end{sidewaystable}

\begin{sidewaystable}[htb]
\tiny
\centering
\setlength\tabcolsep{4pt}
\renewcommand{\arraystretch}{1.2}
\caption{$L_2$ norm of the error for the unsteady advection-diffusion problem}
\label{tab:orderuadf}
\begin{tabular}{llllllllllllllllllllll}
\hline
                   &              & \multicolumn{6}{l}{FR}                                 &  & \multicolumn{6}{l}{HFR}                                &  & \multicolumn{6}{l}{EFR}                                \\ \cline{3-8} \cline{10-15} \cline{17-22} 
$p$                & Mesh size    & $c_{DG}$ & Order & $c_{SD}$ & Order & $c_{HU}$ & Order &  & $c_{DG}$ & Order & $c_{SD}$ & Order & $c_{HU}$ & Order &  & $c_{DG}$ & Order & $c_{SD}$ & Order & $c_{HU}$ & Order \\ \hline
\multirow{4}{*}{1} & 5$\times$5   & 1.72E-02 & -     & 2.28E-02 & -     & 5.19E-02 & -     &  & 1.28E-02 & -     & 2.10E-02 & -     & 5.50E-02 & -     &  & 3.19E-02 & -     & 3.61E-02 & -     & 6.02E-02 & -     \\
                   & 10$\times$10 & 3.68E-03 & 2.23  & 6.10E-03 & 1.90  & 2.02E-02 & 1.36  &  & 2.20E-03 & 2.54  & 5.83E-03 & 1.85  & 2.06E-02 & 1.42  &  & 6.33E-03 & 2.33  & 8.81E-03 & 2.04  & 2.08E-02 & 1.53  \\
                   & 20$\times$20 & 8.33E-04 & 2.14  & 1.52E-03 & 2.00  & 5.62E-03 & 1.84  &  & 4.34E-04 & 2.34  & 1.57E-03 & 1.89  & 5.98E-03 & 1.78  &  & 8.80E-04 & 2.85  & 1.95E-03 & 2.18  & 5.57E-03 & 1.90  \\
                   & 40$\times$40 & 1.90E-04 & 2.13  & 3.89E-04 & 1.97  & 1.49E-03 & 1.91  &  & 9.93E-05 & 2.13  & 4.12E-04 & 1.93  & 1.60E-03 & 1.91  &  & 1.31E-04 & 2.74  & 4.71E-04 & 2.05  & 1.44E-03 & 1.95  \\ \hline
\multirow{4}{*}{2} & 5$\times$5   & 1.02E-03 & -     & 1.42E-03 & -     & 2.36E-03 & -     &  & 7.24E-04 & -     & 1.18E-03 & -     & 2.16E-03 & -     &  & 8.58E-04 & -     & 1.19E-03 & -     & 1.91E-03 & -     \\
                   & 10$\times$10 & 1.24E-04 & 3.03  & 1.64E-04 & 3.12  & 2.47E-04 & 3.26  &  & 8.43E-05 & 3.10  & 1.16E-04 & 3.34  & 1.99E-04 & 3.44  &  & 9.43E-05 & 3.19  & 1.14E-04 & 3.39  & 1.72E-04 & 3.47  \\
                   & 20$\times$20 & 1.52E-05 & 3.04  & 1.85E-05 & 3.15  & 2.56E-05 & 3.27  &  & 1.02E-05 & 3.04  & 1.34E-05 & 3.12  & 2.17E-05 & 3.19  &  & 1.09E-05 & 3.11  & 1.20E-05 & 3.25  & 1.67E-05 & 3.36  \\
                   & 40$\times$40 & 1.90E-06 & 3.00  & 2.22E-06 & 3.06  & 2.91E-06 & 3.14  &  & 1.27E-06 & 3.01  & 1.71E-06 & 2.97  & 2.79E-06 & 2.96  &  & 1.31E-06 & 3.05  & 1.40E-06 & 3.10  & 1.90E-06 & 3.14  \\ \hline
\multirow{4}{*}{3} & 5$\times$5   & 8.04E-05 & -     & 1.22E-04 & -     & 1.64E-04 & -     &  & 5.34E-05 & -     & 7.30E-05 & -     & 1.02E-04 & -     &  & 5.43E-05 & -     & 7.29E-05 & -     & 9.87E-05 & -     \\
                   & 10$\times$10 & 5.21E-06 & 3.95  & 6.73E-06 & 4.18  & 8.39E-06 & 4.29  &  & 3.23E-06 & 4.05  & 4.16E-06 & 4.13  & 5.61E-06 & 4.18  &  & 3.22E-06 & 4.07  & 4.05E-06 & 4.17  & 5.29E-06 & 4.22  \\
                   & 20$\times$20 & 3.22E-07 & 4.01  & 4.14E-07 & 4.02  & 5.12E-07 & 4.03  &  & 2.00E-07 & 4.02  & 2.64E-07 & 3.98  & 3.58E-07 & 3.97  &  & 1.97E-07 & 4.03  & 2.52E-07 & 4.01  & 3.36E-07 & 3.98  \\
                   & 40$\times$40 & 2.07E-08 & 3.96  & 2.53E-08 & 4.03  & 3.05E-08 & 4.07  &  & 1.25E-08 & 4.00  & 1.71E-08 & 3.95  & 2.37E-08 & 3.92  &  & 1.21E-08 & 4.02  & 1.61E-08 & 3.97  & 2.22E-08 & 3.92  \\ \hline
\multirow{4}{*}{4} & 5$\times$5   & 5.37E-06 & -     & 6.98E-06 & -     & 8.28E-06 & -     &  & 3.29E-06 & -     & 4.33E-06 & -     & 5.45E-06 & -     &  & 3.34E-06 & -     & 4.12E-06 & -     & 4.97E-06 & -     \\
                   & 10$\times$10 & 1.64E-07 & 5.04  & 2.21E-07 & 4.98  & 2.65E-07 & 4.97  &  & 1.01E-07 & 5.03  & 1.31E-07 & 5.05  & 1.62E-07 & 5.07  &  & 1.02E-07 & 5.03  & 1.20E-07 & 5.10  & 1.43E-07 & 5.12  \\
                   & 20$\times$20 & 5.16E-09 & 4.99  & 6.75E-09 & 5.04  & 7.96E-09 & 5.05  &  & 3.14E-09 & 5.01  & 4.22E-09 & 4.95  & 5.31E-09 & 4.93  &  & 3.14E-09 & 5.02  & 3.75E-09 & 5.00  & 4.51E-09 & 4.99  \\
                   & 40$\times$40 & 1.64E-10 & 4.98  & 2.07E-10 & 5.03  & 2.40E-10 & 5.05  &  & 9.81E-11 & 5.00  & 1.36E-10 & 4.95  & 1.74E-10 & 4.93  &  & 9.76E-11 & 5.01  & 1.19E-10 & 4.98  & 1.45E-10 & 4.96  \\ \hline
\end{tabular}

\caption{$L_2$-norm of the post-processed solution error for the unsteady linear diffusion problem}
\label{tab:orderzuadf}
\begin{tabular}{lllllllllllllll}
\hline
                   &              & \multicolumn{6}{l}{FR}                                 &  & \multicolumn{6}{l}{HFR}                                \\ \cline{3-8} \cline{10-15} 
$p$                & Mesh size    & $c_{DG}$ & Order & $c_{SD}$ & Order & $c_{HU}$ & Order &  & $c_{DG}$ & Order & $c_{SD}$ & Order & $c_{HU}$ & Order \\ \hline
\multirow{4}{*}{1} & 5$\times$5   & 1.56E-02 & -     & 1.96E-02 & -     & 4.86E-02 & -     &  & 1.14E-02 & -     & 1.85E-02 & -     & 5.08E-02 & -     \\
                   & 10$\times$10 & 3.31E-03 & 2.24  & 5.82E-03 & 1.75  & 1.99E-02 & 1.29  &  & 1.62E-03 & 2.82  & 5.33E-03 & 1.80  & 1.97E-02 & 1.36  \\
                   & 20$\times$20 & 7.11E-04 & 2.22  & 1.45E-03 & 2.01  & 5.58E-03 & 1.83  &  & 2.08E-04 & 2.96  & 1.49E-03 & 1.84  & 5.88E-03 & 1.75  \\
                   & 40$\times$40 & 1.58E-04 & 2.17  & 3.75E-04 & 1.95  & 1.49E-03 & 1.90  &  & 2.62E-05 & 2.99  & 3.96E-04 & 1.91  & 1.58E-03 & 1.89  \\ \hline
\multirow{4}{*}{2} & 5$\times$5   & 9.25E-04 & -     & 1.36E-03 & -     & 2.17E-03 & -     &  & 3.18E-04 & -     & 7.83E-04 & -     & 1.56E-03 & -     \\
                   & 10$\times$10 & 8.51E-05 & 3.44  & 1.23E-04 & 3.46  & 1.87E-04 & 3.53  &  & 2.42E-05 & 3.72  & 5.85E-05 & 3.74  & 1.14E-04 & 3.77  \\
                   & 20$\times$20 & 9.44E-06 & 3.17  & 1.25E-05 & 3.30  & 1.71E-05 & 3.45  &  & 1.66E-06 & 3.86  & 3.97E-06 & 3.88  & 7.68E-06 & 3.90  \\
                   & 40$\times$40 & 1.13E-06 & 3.06  & 1.43E-06 & 3.13  & 1.83E-06 & 3.23  &  & 1.08E-07 & 3.94  & 2.58E-07 & 3.95  & 4.96E-07 & 3.95  \\ \hline
\multirow{4}{*}{3} & 5$\times$5   & 3.70E-05 & -     & 6.32E-05 & -     & 8.51E-05 & -     &  & 1.03E-05 & -     & 2.85E-05 & -     & 4.13E-05 & -     \\
                   & 10$\times$10 & 2.70E-06 & 3.77  & 4.52E-06 & 3.80  & 5.88E-06 & 3.86  &  & 2.46E-07 & 5.39  & 7.15E-07 & 5.32  & 1.00E-06 & 5.36  \\
                   & 20$\times$20 & 1.48E-07 & 4.19  & 2.26E-07 & 4.33  & 2.84E-07 & 4.37  &  & 5.96E-09 & 5.37  & 1.83E-08 & 5.29  & 2.40E-08 & 5.38  \\
                   & 40$\times$40 & 1.10E-08 & 3.75  & 1.47E-08 & 3.94  & 1.75E-08 & 4.02  &  & 1.65E-10 & 5.18  & 5.28E-10 & 5.12  & 6.69E-10 & 5.17  \\ \hline
\multirow{4}{*}{4} & 5$\times$5   & 3.50E-06 & -     & 5.31E-06 & -     & 6.25E-06 & -     &  & 7.69E-07 & -     & 1.58E-06 & -     & 1.83E-06 & -     \\
                   & 10$\times$10 & 8.29E-08 & 5.40  & 1.18E-07 & 5.49  & 1.36E-07 & 5.53  &  & 1.57E-08 & 5.61  & 3.04E-08 & 5.70  & 3.54E-08 & 5.70  \\
                   & 20$\times$20 & 2.52E-09 & 5.04  & 3.37E-09 & 5.13  & 3.79E-09 & 5.16  &  & 2.66E-10 & 5.89  & 5.20E-10 & 5.87  & 6.17E-10 & 5.84  \\
                   & 40$\times$40 & 8.13E-11 & 4.96  & 1.05E-10 & 5.00  & 1.17E-10 & 5.01  &  & 4.26E-12 & 5.97  & 8.41E-12 & 5.95  & 1.01E-11 & 5.93  \\ \hline
\end{tabular}
\end{sidewaystable}
\begin{sidewaystable}[htb]
\tiny
\centering
\setlength\tabcolsep{4pt}
\renewcommand{\arraystretch}{1.2}
\caption{$L_2$ norm of the temperature error for the planar Couette flow case}
\label{tab:ordercouette}
\begin{tabular}{llllllllllllllllllllll}
\cline{3-22}
\multicolumn{1}{l}{} & \multicolumn{1}{l}{} & \multicolumn{6}{l}{FR}                                                           &  & \multicolumn{6}{l}{HFR}                                                           &  & \multicolumn{6}{l}{EFR}                                                                                            \\ \cline{3-8} \cline{10-15} \cline{17-22} 
$p$                  & Mesh size            & $c_{DG}$          & Order & $c_{SD}$          & Order & $c_{HU}$         & Order &  & $c_{DG}$          & Order & $c_{SD}$          & Order & $c_{HU}$          & Order &  & \multicolumn{1}{c}{$c_{DG}$} & Order & \multicolumn{1}{c}{$c_{SD}$} & Order & \multicolumn{1}{c}{$c_{HU}$} & Order \\ \hline
\multirow{4}{*}{1}   & 4$\times$2           & $2.02\text{E-}2$  & -     & $3.02\text{E-}2$  & -     & $7.43\text{E-}2$ & -     &  & $1.87\text{E-}2$  & -     & $1.78\text{E-}2$  & -     & $2.00\text{E-}2$  & -     &  & $1.49\text{E-}2$             & -     & $1.51\text{E-}2$             & -     & $2.09\text{E-}2$             & -     \\
                     & 8$\times$4           & $5.22\text{E-}3$  & 1.95  & $6.69\text{E-}3$  & 2.18  & $1.36\text{E-}2$ & 2.45  &  & $4.82\text{E-}3$  & 1.95  & $4.73\text{E-}3$  & 1.91  & $5.79\text{E-}3$  & 1.79  &  & $3.55\text{E-}3$             & 2.07  & $4.00\text{E-}3$             & 1.92  & $6.65\text{E-}3$             & 1.65  \\
                     & 16$\times$8          & $1.26\text{E-}3$  & 2.05  & $1.47\text{E-}3$  & 2.18  & $2.71\text{E-}3$ & 2.33  &  & $1.10\text{E-}3$  & 2.14  & $1.06\text{E-}3$  & 2.16  & $1.45\text{E-}3$  & 2.00  &  & $9.04\text{E-}4$             & 1.97  & $1.04\text{E-}3$             & 1.94  & $1.74\text{E-}3$             & 1.93  \\
                     & 32$\times$16         & $3.11\text{E-}4$  & 2.02  & $3.59\text{E-}4$  & 2.04  & $6.12\text{E-}4$ & 2.15  &  & $2.42\text{E-}4$  & 2.18  & $2.34\text{E-}4$  & 2.18  & $3.86\text{E-}4$  & 1.91  &  & $2.29\text{E-}4$             & 1.98  & $2.56\text{E-}4$             & 2.03  & $3.96\text{E-}4$             & 2.14  \\ \hline
\multirow{4}{*}{2}   & 4$\times$2           & $1.70\text{E-}3$  & -     & $2.56\text{E-}3$  & -     & $4.01\text{E-}3$ & -     &  & $1.50\text{E-}3$  & -     & $1.55\text{E-}3$  & -     & $1.76\text{E-}3$  & -     &  & $1.53\text{E-}3$             & -     & $1.55\text{E-}3$             & -     & $1.81\text{E-}3$             & -     \\
                     & 8$\times$4           & $2.20\text{E-}4$  & 2.95  & $2.95\text{E-}4$  & 3.12  & $4.33\text{E-}4$ & 3.21  &  & $1.95\text{E-}4$  & 2.95  & $2.07\text{E-}4$  & 2.90  & $2.47\text{E-}4$  & 2.84  &  & $2.00\text{E-}4$             & 2.93  & $2.13\text{E-}4$             & 2.86  & $2.58\text{E-}4$             & 2.81  \\
                     & 16$\times$8          & $2.77\text{E-}5$  & 2.99  & $3.57\text{E-}5$  & 3.05  & $5.01\text{E-}5$ & 3.11  &  & $2.34\text{E-}5$  & 3.06  & $2.54\text{E-}5$  & 3.03  & $3.21\text{E-}5$  & 2.94  &  & $2.46\text{E-}5$             & 3.03  & $2.65\text{E-}5$             & 3.00  & $3.35\text{E-}5$             & 2.95  \\
                     & 32$\times$16         & $3.55\text{E-}6$  & 2.97  & $4.69\text{E-}6$  & 2.93  & $6.56\text{E-}6$ & 2.93  &  & $2.75\text{E-}6$  & 3.09  & $3.09\text{E-}6$  & 3.04  & $4.24\text{E-}6$  & 2.92  &  & $2.92\text{E-}6$             & 3.07  & $3.14\text{E-}6$             & 3.08  & $4.06\text{E-}6$             & 3.04  \\ \hline
\multirow{4}{*}{3}   & 4$\times$2           & $1.55\text{E-}4$  & -     & $2.33\text{E-}4$  & -     & $3.14\text{E-}4$ & -     &  & $1.35\text{E-}4$  & -     & $1.44\text{E-}4$  & -     & $1.61\text{E-}4$  & -     &  & $1.35\text{E-}4$             & -     & $1.45\text{E-}4$             & -     & $1.64\text{E-}4$             & -     \\
                     & 8$\times$4           & $1.08\text{E-}5$  & 3.84  & $1.56\text{E-}5$  & 3.90  & $2.02\text{E-}5$ & 3.96  &  & $9.12\text{E-}6$  & 3.88  & $1.03\text{E-}5$  & 3.80  & $1.21\text{E-}5$  & 3.73  &  & $9.26\text{E-}6$             & 3.87  & $1.05\text{E-}5$             & 3.78  & $1.24\text{E-}5$             & 3.72  \\
                     & 16$\times$8          & $7.13\text{E-}7$  & 3.92  & $9.92\text{E-}7$  & 3.97  & $1.28\text{E-}6$ & 3.98  &  & $5.74\text{E-}7$  & 3.99  & $6.83\text{E-}7$  & 3.92  & $8.39\text{E-}7$  & 3.86  &  & $5.88\text{E-}7$             & 3.98  & $6.93\text{E-}7$             & 3.92  & $8.49\text{E-}7$             & 3.87  \\
                     & 32$\times$16         & $4.62\text{E-}8$  & 3.95  & $6.55\text{E-}8$  & 3.92  & $8.43\text{E-}8$ & 3.92  &  & $3.48\text{E-}8$  & 4.04  & $4.30\text{E-}8$  & 3.99  & $5.50\text{E-}8$  & 3.93  &  & $3.59\text{E-}8$             & 4.04  & $4.26\text{E-}8$             & 4.02  & $5.35\text{E-}8$             & 3.99  \\ \hline
\multirow{4}{*}{4}   & 4$\times$2           & $1.47\text{E-}5$  & -     & $2.21\text{E-}5$  & -     & $2.74\text{E-}5$ & -     &  & $1.29\text{E-}5$  & -     & $1.42\text{E-}5$  & -     & $1.57\text{E-}5$  & -     &  & $1.29\text{E-}5$             & -     & $1.44\text{E-}5$             & -     & $1.60\text{E-}5$             & -     \\
                     & 8$\times$4           & $5.54\text{E-}7$  & 4.73  & $8.24\text{E-}7$  & 4.75  & $1.01\text{E-}6$ & 4.76  &  & $4.75\text{E-}7$  & 4.76  & $5.74\text{E-}7$  & 4.63  & $6.61\text{E-}7$  & 4.57  &  & $4.79\text{E-}7$             & 4.75  & $5.82\text{E-}7$             & 4.62  & $6.74\text{E-}7$             & 4.57  \\
                     & 16$\times$8          & $1.89\text{E-}8$  & 4.88  & $2.82\text{E-}8$  & 4.87  & $3.46\text{E-}8$ & 4.87  &  & $1.55\text{E-}8$  & 4.94  & $1.98\text{E-}8$  & 4.86  & $2.35\text{E-}8$  & 4.81  &  & $1.57\text{E-}8$             & 4.93  & $2.02\text{E-}8$             & 4.85  & $2.41\text{E-}8$             & 4.81  \\
                     & 32$\times$16         & $6.17\text{E-}10$ & 4.93  & $9.34\text{E-}10$ & 4.91  & $1.15\text{E-}9$ & 4.91  &  & $4.75\text{E-}10$ & 5.03  & $6.22\text{E-}10$ & 4.99  & $7.58\text{E-}10$ & 4.95  &  & $4.86\text{E-}10$            & 5.02  & $6.34\text{E-}10$            & 4.99  & $7.70\text{E-}10$            & 4.97  \\ \hline
\end{tabular}
\end{sidewaystable}
\FloatBarrier

\end{document}